\documentclass[preprint, 11pt]{elsarticle}

\usepackage[T1]{fontenc}
\usepackage[utf8]{inputenc}
\usepackage[ngerman,english]{babel}
\usepackage[hmargin=0.9in]{geometry}
\usepackage[babel,german=quotes]{csquotes}
\usepackage{subfig}
\usepackage{graphicx}
\usepackage[colorlinks=true,citecolor=blue,urlcolor=blue]{hyperref}
\usepackage{caption}
\usepackage{amsmath}          
\usepackage{amsthm}           
\usepackage{amssymb}          
\usepackage{bm, amsfonts}     
\usepackage{xcolor}
\usepackage{fixmath}
\usepackage{tikz}
\usepackage{bm}
\usepackage{makecell}

\newcommand {\dx} {\,{\rm d}{\mathbf x}}

\newcommand {\ds} {\,{\rm d}{\mathrm s}}

\newcommand {\bu} {{\mathbf v}}

\newcommand {\bc} {{\mathbf c}}

  \newcommand {\bF} {{\mathbf f}}

   \newcommand{\bff}{\boldsymbol{f}}

  \newcommand{\R}{\mathbb{R}}
  \newdefinition{rmk}{Remark}
  
  \newcommand{\pd}[2]{\frac{\partial #1}{\partial #2}}

  \newcommand{\td}[2]{\frac{\mathrm d #1}{\mathrm d #2}}

\newcommand{\beq}{\begin{equation}}
\newcommand{\eeq}{\end{equation}}

\newcommand\red[1]{{\color{black}#1}} 
\newcommand\cyan[1]{{\color{black}#1}}  
\newcommand\black[1]{{\color{black}#1}}	

\newcommand{\tc}{\tilde{\mathbf{c}}}
\newcommand{\Q}{\mathbb{Q}_1}
\newcommand{\QQ}{\mathbb{Q}_2}
\newcommand{\bfF}{\mathbf{F}}
\newcommand{\bfx}{{\bf x}}

\makeatletter
\def\ps@pprintTitle{%
  \let\@oddhead\@empty
  \let\@evenhead\@empty
  \def\@oddfoot{
    \footnotesize\itshape
    \hfill\today
  }%
  \let\@evenfoot\@oddfoot}
\makeatother

\begin{document}

\begin{frontmatter} 
  \title{Subcell flux limiting for high-order Bernstein finite element discretizations of scalar hyperbolic conservation laws}

\author[TUDo]{Dmitri Kuzmin\corref{cor1}}
\ead{kuzmin@math.uni-dortmund.de}
\cortext[cor1]{Corresponding author}

\address[TUDo]{Institute of Applied Mathematics (LS III), TU Dortmund University\\ Vogelpothsweg 87,
  D-44227 Dortmund, Germany}

\author[KAUST]{Manuel Quezada de Luna}
\ead{manuel.quezada@kaust.edu.sa}

\address[KAUST]{King Abdullah University of Science and Technology (KAUST)\\ Thuwal 23955-6900, Saudi Arabia}

\journal{Computers and Mathematics with Applications}

\begin{abstract}
  This \black{work} extends the concepts of algebraic flux correction and convex limiting
  to continuous high-order Bernstein finite element discretizations of scalar
  hyperbolic problems. Using an array of adjustable diffusive fluxes, the
  standard Galerkin approximation is transformed into a nonlinear high-resolution
  scheme which has the compact sparsity pattern of  the piecewise-linear or
  multilinear subcell discretization. The representation of this scheme in
  terms of invariant domain preserving states makes it possible to
  prove the validity of local discrete maximum principles under CFL-like
  conditions. In contrast to predictor-corrector approaches based on
  the flux-corrected transport methodology, the proposed flux limiting strategy
  is monolithic, i.e., limited antidiffusive terms are incorporated into the
  well-defined residual of a nonlinear (semi-)discrete problem. A stabilized high-order
  Galerkin discretization is recovered
  if no limiting is performed. In the limited version, the
  compact stencil property \black{prevents} direct mass exchange between nodes that
  are not nearest neighbors. A formal proof of sparsity is provided for
  simplicial and box elements. The involved element contributions 
  can be calculated efficiently
  making use of matrix-free algorithms and precomputed element matrices
  of the reference element. Numerical studies for
  \black{$\mathbb{Q}_2$} discretizations of \black{linear and nonlinear}
  two-dimensional test problems illustrate the virtues of monolithic convex
  limiting based on subcell flux decompositions.
  
\end{abstract}
\end{frontmatter}
\pagebreak

\section{Introduction}

Algebraic flux correction (AFC) \cite{afc_analysis1,afc_analysis2,afc1,CL-diss} is a general
framework for the design of bound-preserving finite element schemes. Many representatives
of nonlinear high-resolution AFC schemes are based on algebraic interpretations and
generalizations of flux-based structured grid methods for hyperbolic
conservation laws. Finite element AFC versions of upwinding techniques, flux-corrected
transport (FCT) algorithms \cite{fct1,zalesak79}, total variation diminishing (TVD)
limiters \cite{harten1,harten2}, and their local extremum diminishing (LED) counterparts
\cite{jameson1,jameson2} have been used since the late 1980s \cite{badia,afc1,fctools,lohner87,baum1994,tvd,peraire1993,Selmin1993,Selmin1996}. In recent years, their further
development was stimulated by major breakthroughs in theoretical analysis of
the involved `variational crimes'. The work of Barrenechea et al. \cite{Barrenechea2016,afc_analysis1,afc_analysis2} established a theoretical framework for proving convergence and well-posedness of AFC schemes for steady convection-diffusion equations. Lohmann \cite{CL-diss} extended this framework to finite element discretizations of steady and unsteady linear advection
problems.  Guermond et al.~\cite{Guermond2014,Guermond2018,Guermond2016,Guermond2017} introduced a family of explicit invariant domain preserving (IDP) schemes for nonlinear hyperbolic problems. Their analytical studies paved the way for the development of novel convex limiting techniques \cite{Guermond2018,Guermond2019,convex} based on generalizations of localized FCT schemes \cite{cotter,CG-BFCT} and monolithic AFC approaches \cite{convex}.

As of this writing, the overwhelming majority of algebraic flux correction tools and
the underlying theory are not readily applicable to finite element approximations of
degree $p>1$. Using the Bernstein basis representation, a few element-based high-order
extensions of residual distribution methods \cite{trefilik,RD-BFCT} and localized FCT
schemes \cite{DG-BFCT,CG-BFCT} were developed for continuous and discontinuous Galerkin
discretizations. A common drawback of the underlying limiting techniques for
antidiffusive element contributions is the possibility of direct mass exchange between
all nodes of a high-order Bernstein element. This lack of locality was
found to be acceptable in applications to linear advection problems
\cite{DG-BFCT,RD-BFCT,CG-BFCT} but the
design of high-resolution AFC schemes
for nonlinear conservation laws calls for the use of
flux-based subcell approximations with compact computational stencils.

\black{
The AFC methodology that we introduce in the present paper converts a high-order
continuous Galerkin discretization into a nonlinear IDP scheme with the compact
sparsity pattern of a piecewise $\mathbb{P}_1/\mathbb{Q}_1$ subcell approximation. We begin in \S\ref{sec:high} with the description of the high-order
  Bernstein finite element discretization.} \black{Then, in \S\ref{sec:low},
  we derive a low-order IDP approximation which has a compact stencil and
  is less diffusive than the full stencil version using the same kind of
  algebraic residual correction (discrete upwinding \cite{afc1,fctools,CG-BFCT}
  or Rusanov dissipation \cite{Guermond2016,RD-BFCT,convex,afc2}).
Next, in \S\ref{sec:convex}, we present a monolithic convex
  limiting procedure for the antidiffusive correction
  terms corresponding to a (stabilized) high-order target. The compact
  stencil property is preserved using a decomposition of the antidiffusive
  element
contributions into subcell fluxes between nearest neighbor nodes.
This approach, which is described in \S\ref{sec:subcell}, involves
the solution of small sparse linear systems on each macroelement. The IDP
  property of the corresponding discrete problem is shown using the
  proof techniques developed in \cite{Guermond2018,convex}.
  In \S\ref{sec:stabilization} and \S\ref{sec:smoothness_indicator}, we discuss the optional stabilization techniques for the high-order target flux and
  \red{Laplacian}-based smoothness indicators that preserve the high-order accuracy near smooth local extrema.
  Time integration is performed using an explicit
  (third order with three stages)
    strong stability preserving Runge-Kutta method \cite{ssprev, ssprev0}.
    The possibility of using precomputed element matrices of the reference element and
matrix-free solvers for the global system may be exploited in efficient implementations
of the proposed algorithms.
  The results of numerical studies for linear and nonlinear conservation
  laws are presented in \S\ref{sec:num}.
  Finally, we close in \S\ref{sec:conclusions} with conclusions.}

\section{High-order Bernstein finite element discretization}
\label{sec:high}

\black{We} restrict our presentation to the case of a  scalar 
conservation law. An extension of the proposed methodology to nonlinear
hyperbolic systems can be carried out as in \cite{convex}
\black{and will be presented elsewhere}. 
Let $u(\mathbf{x},t)$ be a scalar
quantity of interest depending on the space location $\mathbf{x}\in
\R^d,\ d\in\{1,2,3\}$ and time instant $t\ge 0$. Consider
an initial-boundary value problem of the form \cite{Guermond2018,convex}
\begin{subequations}
\begin{align}
 \pd{u}{t}+\nabla\cdot\mathbf{f}(u)=0 &\qquad\mbox{in}\ \Omega\times\R_+,
\label{ibvp-pde}\\
 u(\cdot,0)=u_0 &\qquad\mbox{in}\ \Omega,\label{ibvp-ic}\\
 (u-u_{\rm in})\mathbf{f}'(u)\cdot\mathbf{n}=0
  &\qquad\mbox{on}\ \Gamma_-,
\label{ibvp-bc}
\end{align}
\end{subequations}
where $\Omega\subset\R^d$ is a bounded domain, 
$\mathbf{f}=(\mathsf{f}_1,\ldots,\mathsf{f}_d)$ is a possibly
nonlinear flux function,
$u_0$ is the initial data, $u_{\rm in}$ is the Dirichlet
boundary data, $\mathbf{n}$ is the unit outward
normal to the Lipschitz boundary $\Gamma=\partial\Omega$, and
$\Gamma_-=\{\mathbf{x}\in\Gamma\,:\,\mathbf{f}'(u)\cdot
\mathbf{n}<0\}$ is the hyperbolic inlet.

Suppose that the exact solution $u$ belongs to a convex set \cyan{$\mathcal G\subset\R$} for all $t\ge 0$. Then $\mathcal G$ is called an invariant set of problem  \eqref{ibvp-pde}--\eqref{ibvp-bc}, and it is natural to require that numerical approximations belong to (a subset of)  $\mathcal G$ as well. Adopting the terminology of Guermond et al. \cite{Guermond2018,Guermond2016,Guermond2017}, we will call 
a discretization of problem \eqref{ibvp-pde}--\eqref{ibvp-bc} invariant
domain preserving (IDP) if the solution of the (semi-) discrete problem
is guaranteed to stay in a convex invariant set.

To begin with, we discretize \eqref{ibvp-pde} in space using a high-order
continuous Galerkin method.
Given a conforming mesh $\mathcal T_h=\{K^1,\ldots,K^{E_h}\}$, we
define a finite element
approximation $u_h\approx u$ in terms of globally continuous piecewise-polynomial
basis functions $\varphi_j$, where $j\in\{1,\ldots, N_h\}$ is the global number of a
nodal point $\mathbf{x}_j$. The local number $j_e=\mathcal I^e(j)$ of node $j$
in $K^e$ is determined by a mapping $\mathcal I^e:\{1,\ldots,N_h\}\to\{1,\ldots,N\}.$
The corresponding local basis function is denoted by $\varphi_{j_e}^e$.
The global numbers of nodes $\mathbf{x}_1^e,\ldots,\mathbf{x}_N^e$
belonging to $K^e$ are stored in the integer set
$\mathcal N^e\subset\{1,\ldots, N_h\}$. 

The polynomial restriction of $u_h=\sum_{j=1}^{N_h}u_j\varphi_j$ to element $K^e,\ e=1,\ldots,E_h$ is given by
\beq\label{uh-def} 
u_h^e:=u_h|_{K^e}
=\sum_{j\in\mathcal N^e}u_j\varphi_j
=\sum_{j\in\mathcal N^e}u_{j}\varphi_{j_e}^e
=\sum_{i=1}^Nu_i^e\varphi_i^e,
\eeq
where $u_i^e=u_{j_e}$ is the degree of freedom (DoF) associated
with the nodal point $\mathbf{x}_i^e=\mathbf{x}_{j_e},\ j\in\mathcal N^e$. 

To enforce the IDP property using algebraic flux correction
\cite{afc_analysis1,afc_analysis2,afc1,CL-diss} in what follows,
we will use the Bernstein basis representation of $u_h$.
The Bernstein basis functions $\varphi_j^e$, the definition of
which for simplicial and tensor product meshes can be found in
the Appendix, are nonnegative and form a partition of unity,
i.e., $\sum_{j=1}^N\varphi_j^e\equiv 1$. It follows that for
any $\mathbf{x}\in K^e$, the state
$u_h(\mathbf{x})$ is a convex combination of the
nodal states $u_1^e,\ldots,u_N^e$. Thus, we have
\beq
u_1^e,\ldots,u_N^e\in\mathcal G\quad\Rightarrow\quad u_h(\mathbf{x})\in
\mathcal G\quad \forall\mathbf{x}\in K^e
\eeq
for any convex invariant set $\mathcal G$ of the hyperbolic initial-boundary
value problem \eqref{ibvp-pde}--\eqref{ibvp-bc}.\smallskip

\cyan{Integrating the weighted residuals of \eqref{ibvp-pde} and
  \eqref{ibvp-bc} over $\Omega$ and $\Gamma_-$, respectively, we
  obtain a weak form of the problem at hand. The standard 
Galerkin discretization replaces it with 
\beq\label{weak1}
\sum_{e=1}^{E_h}
\int_\Omega w_h \left(\pd{u_h}{t}+\nabla\cdot\mathbf{f}(u_h)\right)\dx
=\sum_{e=1}^{E_h}
\int_{\partial K^e\cap \Gamma_-} w_h(u_h-u_{\rm in})\mathbf{f}'(u_h)\cdot\mathbf{n}
\ds\qquad \forall w_h\in W_h,
\eeq
where $W_h$ is the finite-dimensional space spanned by the
Bernsiein basis functions $\varphi_1,\ldots,\varphi_{N_h}$.
}

Substitution of \eqref{uh-def} into \eqref{weak1} with the
test function $w_h=\varphi_i$
produces the semi-discrete equation
\beq
\label{galerkin1}
\sum_{j\in\mathcal N_i}m_{ij}\td{u_j}{t}=b_i(u_h,u_{\rm in})
-\sum_{e\in\mathcal E_i}\int_{K^e}
\varphi_i\nabla\cdot\mathbf{f}(u_h)\dx
,
\eeq
where $\mathcal E_i$ is the set of elements containing node $i$
 and $\mathcal N_i$ is the set of nodes belonging
to these elements. The entries $m_{ij}$ of the global consistent mass matrix
and the boundary term $b_i$ are defined by
\beq\label{mass}
m_{ij}=\sum_{e\in\mathcal E_i\cap\mathcal E_j}m_{ij}^e,\qquad
m_{ij}^e=\int_{K^e}\varphi_i\varphi_j\dx,
\eeq
\beq
b_i(u_h,u_{\rm in})
=\sum_{e\in\mathcal E_i}\int_{\partial K^e\cap\Gamma_-}\varphi_i
(u_h-u_{\rm in})\mathbf{f}'(u_h)\cdot\mathbf{n}\ds.
\eeq
In \black{practice}, only the $N^2$ nonvanishing entries of element matrices
like $M^e_C=\{m_{ij}^e\}_{i,j=1}^{N_h}$ are calculated and inserted into global matrices. To
avoid conversion between global and local indices, we will use the global index
notation for element matrices and vectors in this paper.


\section{Low-order Bernstein finite element discretization}
\label{sec:low}

\begin{figure}[h]
  \centering
  \subfloat[$\mathbb{P}_2$ and $\mathbb{P}_3$ elements]
           {\includegraphics[scale=0.3]{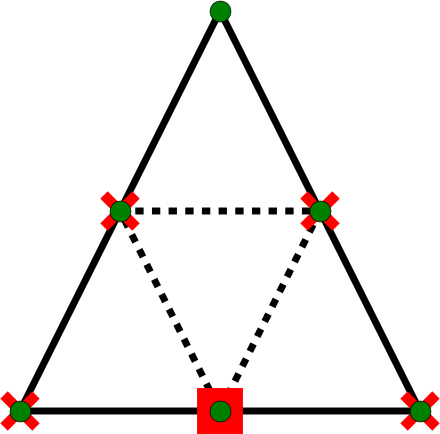}\qquad
            \includegraphics[scale=0.3]{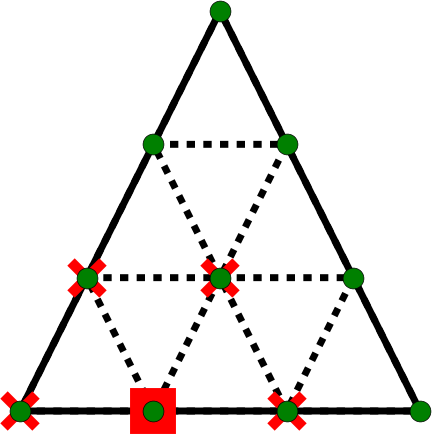}}
  \hspace{40pt}
  \subfloat[$\QQ$ and $\mathbb{Q}_3$ elements]
           {\includegraphics[scale=0.3]{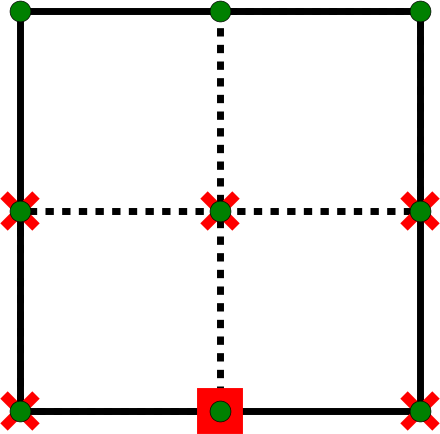}\qquad
             \includegraphics[scale=0.3]{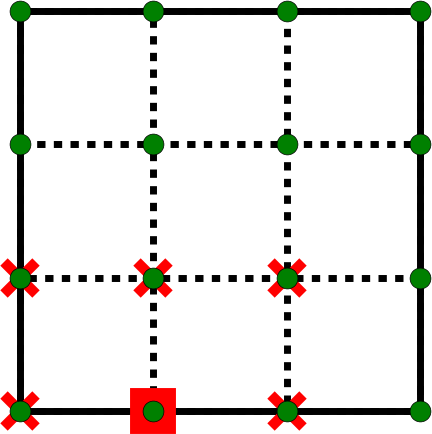}}
           \caption{Nodes and
             subcells of typical high-order elements. The boundary of
             the
             macroelement $K^e$ is marked with solid black lines. The internal
             boundaries of its subcells are marked with dashed black lines.
             All local DoFs are marked with green circles. The red crosses
             correspond to the nearest neighbors of the DoF marked by the
             red square.
      \label{fig:subcells}}
\end{figure}

A space discretization of the form \eqref{galerkin1} can be transformed into
a compact-stencil IDP scheme by using row-sum mass lumping and modifying
the Galerkin element contributions
\beq\label{rhs_ele}
\int_{\partial K^e\cap\Gamma_-}\varphi_i
(u_h-u_{\rm in})\mathbf{f}'(u_h)\cdot\mathbf{n}\ds
-\int_{K^e}\varphi_i\nabla\cdot\mathbf{f}(u_h)\dx.
\eeq
Approximating the flux $\mathbf{f}(u_h)$ by the {\it group finite element} interpolant
\cite{group,fletcher1,fletcher2,Selmin1993,Selmin1996}
 \beq\label{groupf}
  \mathbf{f}_h^e=\sum_{j\in\mathcal N^e}\mathbf{f}_j\varphi_j,
  \qquad \mathbf{f}_j=(\mathsf{f}_{j,1},\ldots,\mathsf{f}_{j,d})=
  \mathbf{f}(u_j)
  \eeq
  and using a lumped approximation of the boundary term,
  we replace \eqref{rhs_ele} with
  \beq
  \int_{\partial K^e\cap\Gamma_-}
   \varphi_i(u_i-u_{\rm in})
   \mathbf{f}'(u_h)\cdot\mathbf{n}\ds
   -\sum_{j\in\mathcal N^e}\mathbf{c}_{ij}^e\cdot\mathbf{f}_j.
   \eeq
 The vector valued coefficients
 $\bc_{ij}^e=(c_{ij,1}^e,\ldots,c_{ij,d}^e)$ of the discrete gradient
operator are defined by 
\beq\label{cij-def}
\bc_{ij}^e=\int_{K^e}\varphi_i\nabla\varphi_j\dx=-\mathbf{c}_{ji}^e
+\int_{\partial K^e}\varphi_i\varphi_j\mathbf{n}\ds.
\eeq
The transformation of the consistent element mass matrix $M^e_C=
\{m_{ij}^e\}_{i,j=1}^{N_h}$ into its lumped counterpart $M^e_L=
\{\delta_{ij}m_i^e\}_{i,j=1}^{N_h}$ with the
diagonal entries
\beq
m_i^e=\sum_{j=1}^{N_h}m_{ij}^e=\sum_{j\in\mathcal N_i}m_{ij}^e
=\int_{K^e}\varphi_i\dx=\frac{|K^e|}{N}>0
\eeq
corresponds to multiplication by the local mass lumping
operator
\beq
P^e=M^e_L(M^e_C)^{-1}.\label{prec}
\eeq
Following the approach proposed in \cite{CG-BFCT},
we apply $P^e$ to $C^e_k=\{c_{ij,k}^e\}_{i,j=1}^{N_h},\ k=1,\ldots,d$
as well. As shown in \cite{CG-BFCT} for the 1D case, this
modification produces {\bf sparse} element matrices
\beq\label{precC}
\tilde C_k^e=P^eC_k^e,\qquad k=1,\ldots,d
\eeq
such that $c_{ij,k}^e=0$ for $j\notin\tilde{\mathcal N}_i^e$, where
$\tilde{\mathcal N}_i^e\subseteq\mathcal N^e$ is the local stencil of the $\mathbb{P}_1/\mathbb{Q}_1$ subcell discretization (see Fig.~\ref{fig:subcells}), i.e., the integer set containing
the local numbers of the nearest neighbors of node $i$ in $K^e$. In the
Appendix, we show the compact-stencil property of the element contributions
 $\tilde{\mathbf{C}}^e
=(\tilde C_{1}^e,\ldots,\tilde C_{d}^e)$ to the
lumped discrete gradient operator
for $d$-simplex and $d$-box Bernstein elements.

\begin{rmk}
  Strict positivity of all lumped mass matrix entries $m_i$ and the compact
  sparsity pattern of $\tilde{\mathbf{C}}^e$ are due to the use of the
  Bernstein basis. High-order Lagrange finite elements do not
  provide these properties which will play an important role in the
  derivation of the proposed correction procedures.
\end{rmk}

\begin{rmk}
  In \cite{CG-BFCT} and \cite{RD-BFCT}, the mass lumping operator $P^e$ was applied
  to the element matrix of the advective term discretized without using the
  group finite element formulation \eqref{groupf} for the linear flux function
  $\bff(\mathbf{x},u)=\bu(\mathbf{x}) u$. This approach does not guarantee exact sparsity
  for general velocity fields $\bu(\mathbf{x})$. As a consequence, the resulting
  schemes become less accurate as the polynomial degree $p$ is increased while
  keeping the total number of DoFs $N_h$ fixed \cite{RD-BFCT}.
\end{rmk}

The replacement of $M^e_C$ and $\mathbf{C}^e=(C_{1}^e,\ldots,C_{d}^e)$ with
the lumped element matrices $M^e_L$ and  $\tilde{\mathbf{C}}^e$ is not
enough to guarantee that the modified Galerkin scheme is IDP. To enforce
the IDP property in a provable manner, we replace the element vector
$\tilde{\mathbf{C}}^e\cdot\mathbf{f}^e=\sum_{k=1}^dC^e_k\mathsf{f}^e_k$ 
by $\tilde{\mathbf{C}}^e\cdot\mathbf{f}^e-\tilde{D}^eu^e$, where $\tilde{D}^e=
\{\tilde d_{ij}^e\}_{i,j=1}^{N_h}$ is the element matrix of a graph Laplacian
(discrete diffusion) operator.
\medskip

The above manipulations convert \eqref{galerkin1} into
the compact-stencil low-order approximation
\beq
  m_i\td{u_i}{t}=\sum_{e\in\mathcal E_i}\left(
  \sum_{j\in\tilde{\mathcal N}_i^e\backslash\{i\}}
  \tilde d_{ij}^e(u_j-u_i)- \sum_{j\in\tilde{\mathcal N}_i^e}
  \tilde{\mathbf{c}}_{ij}^e\cdot\mathbf{f}_j\right)
  +\tilde b_i(u_h,u_{\rm in}),
\label{nobar} 
\eeq
where $m_i=\sum_{e\in\mathcal E_i}m_i^e$ is a diagonal entry of the
global lumped mass matrix and
\beq
   \tilde b_i(u_h,u_{\rm in})=
   \sum_{e\in\mathcal E_i}\int_{\partial K^e\cap\Gamma_-}
   \varphi_i(u_i-u_{\rm in})
   \mathbf{f}'(u_h)\cdot\mathbf{n}\ds.
   \eeq

   To define artificial diffusion coefficients \black{$\tilde{d}_{ij}^e$} that guarantee
   the IDP property for general hyperbolic problems, we
   write  \eqref{nobar} in the equivalent form
    \beq  m_i\td{u_i}{t}=\sum_{e\in\mathcal E_i}
 \sum_{j\in\tilde{\mathcal N}_i^e\backslash\{i\}}
2\tilde d_{ij}^e(\bar u_{ij}^e-u_i)+\tilde b_i(u_h,u_{\rm in}),\label{bar}
\eeq
where
   \beq
 \bar u_{ij}^e=\frac{u_i+u_j}{2}-\frac{\tilde{\mathbf{c}}_{ij}^e
   \cdot(\mathbf{f}_j-\mathbf{f}_i)}{2\tilde d_{ij}^e}.
 \eeq
Guermond and Popov \cite{Guermond2016} were the first to recognize that 
representations of explicit schemes in terms of the bar states
$\bar u_{ij}$ lead to remarkably simple proofs of the IDP property. Indeed,
\eqref{bar} exhibits the structure of a discretized diffusion
equation in which the nodal state $u_j\in\mathcal G$ is replaced
with $\bar u_{ij}^e\in\mathcal G$.

If time
discretization is performed using an explicit SSP Runge-Kutta
method \cite{ssprev}, each stage is a forward Euler update of the form
 \beq\label{feuler}
 m_i\bar u_i=m_iu_i+\Delta t\sum_{e\in\mathcal E_i}
 \sum_{j\in\tilde{\mathcal N}^e_i\backslash\{i\}}
2\tilde d_{ij}^e(\bar u_{ij}-u_i)+\tilde b_i(u_h,u_{\rm in}).
\eeq
\red{The result is IDP} 
for time steps $\Delta t$ satisfying the CFL-like condition
\beq\label{cfl}
\Delta t\left(\sum_{e\in\mathcal E_i}
\sum_{j\in\tilde{\mathcal N}^e_i\backslash\{i\}}2\tilde d_{ij}^e
-\sum_{e\in\mathcal E_i}\int_{\partial K^e\cap\Gamma_-}
\varphi_i\mathbf{f}'(u_h)\cdot\mathbf{n}\ds
\right)
\le m_i
\eeq
provided that all $\bar u_{ij}^e$ stay in $\mathcal G$
for $u_i, u_j\in \mathcal G$. As explained
in  \cite{Guermond2016}, this requirement can be satisfied by using the
guaranteed maximum speed (GMS)
  \beq
  \lambda_{ij}^e=
  \max_{\omega\in[0,1]}\left|\mathbf{n}_{ij}^e\cdot
  \mathbf{f}'(\omega u_i+(1-\omega)u_j)\right|,\qquad
\black{\tilde{\mathbf{n}}_{ij}^e}=
  \frac{\tilde{\mathbf{c}}_{ij}^e}{|\tilde{\mathbf{c}}_{ij}^e|}\label{gmsdef}
 \eeq
to define the \black{Rusanov-type} artificial viscosity coefficients
 \beq\label{dij-def}
 \tilde{d}_{ij}^e=\begin{cases}
  \max\{|\tilde{\bc}_{ij}^e|,|\tilde{\bc}_{ji}^e|\}
  \max\{\lambda_{ij}^e,\lambda_{ji}^e\}
& \mbox{if}\ \cyan{i\in\mathcal N^e,\ j\in\mathcal N^e\backslash\{i\},}\\
  -\sum_{k\in \tilde{\mathcal N}_i^e\backslash\{i\}}  \tilde{d}_{ik}^e
  & \mbox{if}\ \cyan{j=i\in\mathcal N^e,}\\
  \cyan{0} & \cyan{\mbox{otherwise}}
  \end{cases}
  \eeq
  such that \cite{convex}
  \beq
  \min\{u_i,u_j\}\le \bar u_{ij}^e
\le \max\{u_i,u_j\}.\label{maxbarstate}
\eeq
Note that the element matrix $\tilde{D}^e$ has the same compact sparsity pattern as
$\tilde{\mathbf{C}}^e$.
\medskip

\black{For linear flux functions of the form $\bff(\mathbf{x},u)=\bu(\mathbf{x}) u$, where
  $\mathbf{v}$ is a spatially variable velocity field, 
  the validity of \eqref{maxbarstate} cannot be guaranteed,
  e.g., in the case when $u_i=u_j$ and
  $\mathbf{v}_i\ne\mathbf{v}_j$ \cite{convex}. \red{The} edge contributions of
  the low-order scheme defined by \eqref{feuler} and
  \eqref{dij-def} are given by
    \begin{align*}
      2 \tilde d_{ij}^e(\bar u_{ij}^e-u_i)
    &= \tilde d_{ij}^e
     (u_j-u_i)-\tilde{\mathbf{c}}_{ij}^e
      \cdot(\mathbf{v}_ju_j-\mathbf{v}_iu_i)\\
      &=\underbrace{(\tilde d_{ij}^e
      -\tilde{\mathbf{c}}_{ij}^e\cdot\mathbf{v}_j)}_{\in[0,2\tilde d_{ij}^e]}u_j
      -\underbrace{(\tilde d_{ij}^e-\tilde{\mathbf{c}}_{ji}^e\cdot\mathbf{v}_i)
       }_{\in[0,2\tilde d_{ij}^e]}u_i.
      \end{align*}
\red{Adapting the
    GMS formula \eqref{gmsdef} to the case of linear advection,
    the maximum speeds that appear in definition \eqref{dij-def}
    of the Rusanov diffusion coefficient $d_{ij}^e$ can be redefined as
    $
\lambda_{ij}^e=
  \max_{\mathbf{x}\in K^e}|\mathbf{v}(\mathbf{x})|$. The resulting}
 approximation is IDP w.r.t.
$\mathcal G=\{u\in\R\,|\, u\ge 0\}$
under the time step restriction \eqref{cfl}.  

  A less dissipative low-order scheme for the 
 linear advection equation  can be
  constructed using
      \begin{align}\label{dij-linear-def}
      \tilde{d}^e_{ij} =
      \begin{cases}
        \max\{\tc^e_{ij}\cdot \bu_j, ~0, ~\tc^e_{ji}\cdot\bu_i\} & \mbox{if }\
       \cyan{i\in\mathcal N^e,\ j\in\mathcal N^e\backslash\{i\},} \\
        -\sum_{k\in \tilde{\mathcal N}_i^e\backslash\{i\}}  \tilde{d}_{ik}^e & \mbox{if }\ \cyan{j=i\in\mathcal N^e,}\\
  \cyan{0} & \cyan{\mbox{otherwise}.}
      \end{cases}
      \end{align}
      This alternative to \eqref{dij-def}
       is known as {\it discrete upwinding} \cite{afc1,fctools,CG-BFCT}. In view of the fact that
      \begin{align*}
      2 \tilde d_{ij}^e(\bar u_{ij}^e-u_i)
    & =\max\{\tc^e_{ij}\cdot \bu_j, ~0, ~\tc^e_{ji}\cdot\bu_i\}
     (u_j-u_i)-\tilde{\mathbf{c}}_{ij}^e
      \cdot(\mathbf{v}_ju_j-\mathbf{v}_iu_i)\\
      &=\underbrace{(\max\{\tc^e_{ij}\cdot \bu_j, ~0, ~\tc^e_{ji}\cdot\bu_i\}
      -\tilde{\mathbf{c}}_{ij}^e\cdot\mathbf{v}_j)}_{\ge 0}u_j\\
      &\phantom{=} -( \max\{\tc^e_{ij}\cdot \bu_j, ~0, ~\tc^e_{ji}\cdot\bu_i\}-
      \tilde{\mathbf{c}}_{ij}^e\cdot\mathbf{v}_i)u_i,
      \end{align*}
      the corresponding low-order scheme \eqref{feuler} is positivity-preserving for sufficiently small time steps $\Delta t$.  \red{It is at most as diffusive as the one based on \eqref{dij-def} since $|\tc^e_{ij}\cdot \bu_j|\le|\tc^e_{ij}|\max_{\mathbf{x}\in K^e}|\mathbf{v}(\mathbf{x})|=|\tc^e_{ij}|\lambda_{ij}^e$.}
      \smallskip
      
      In  \S\ref{sec:num}, we solve linear advection problems using \eqref{dij-linear-def}. For nonlinear conservation laws, we use the GMS formula \eqref{dij-def}. As remarked by Guermond and Popov \cite{Guermond2016}, the use of \eqref{dij-linear-def} with the nodal speeds $\mathbf{v}_i:=\mathbf{f}'(u_i)$ may result in entropy-violating weak solutions to nonlinear problems. 
      
\begin{rmk}
  Instead of assembling the global graph Laplacian $\tilde D$ from sparse element matrices $\tilde{D}^e$ defined
  by \eqref{dij-def} or \eqref{dij-linear-def},
  the global discrete gradient \black{operator} $\tilde{\mathbf{C}}$
  can be used to generate $\tilde D$ after the element-by-element assembly
  from  $\tilde{\mathbf{C}}^e$, cf.  \cite{Guermond2018,convex}. 
\end{rmk}}

\begin{rmk}
The use of explicit SSP Runge-Kutta time discretizations is not a
necessary condition for provable preservation of invariant domains.
However, the verification of IDP properties for implicit
and stationary versions of our low-order scheme requires more
sophisticated analysis (cf. \cite{afc_analysis1,afc_analysis2,CL-diss}).
\end{rmk}

 As we show in the next section,
the bar state form \eqref{bar} of \eqref{nobar} is also ideally suited for
the derivation of high-order extensions that preserve the IDP property using
built-in flux limiters.

 \section{Convex limiting for high-order subcell fluxes}\label{sec:convex}

 Decomposing \eqref{galerkin1} into the
 low-order IDP part \eqref{nobar} and a remainder, we write
 it in the form 
 \beq
  m_i\td{u_i}{t}=\sum_{e\in\mathcal E_i}\left(\sum_{j\in\tilde{\mathcal N}_i^e\backslash\{i\}}
  \black{\tilde{d}_{ij}^e}(u_j-u_i)- \sum_{j\in\tilde{\mathcal N}_i}
  \tilde{\mathbf{c}}_{ij}^e\cdot\mathbf{f}_j+f_i^e+g_i^e\right)
  +\tilde b_i(u_h,u_{\rm in}),
\eeq
where
 \begin{align}
f_i^e&=\textcolor{black}{
\sum_{j\in\tilde{\mathcal N}_i^e\backslash\{i\}} \black{\tilde{d}_{ij}^e}(u_i-u_j)+
\sum_{j\in\mathcal N^e\backslash\{i\}} m_{ij}^e(\dot u_i-\dot u_j)}
+\sum_{\textcolor{black}{j\in\mathcal N^e}}
(\tilde{\mathbf{c}}_{ij}^e-\mathbf{c}_{ij}^e)\cdot\mathbf{f}_j
-\sum_{\textcolor{black}{j\in\mathcal N^e}}
\mathbf{c}_{ji}^e\cdot\mathbf{f}_j
\nonumber\\
&\qquad\qquad+\int_{K^e}\nabla\varphi_i\cdot\mathbf{f}(u_h)\dx
\textcolor{black}{\ +
\int_{\partial K^e\cap\Gamma_-}\varphi_i(u_h-u_i)
   \mathbf{f}'(u_h)\cdot\mathbf{n}\ds},
  \label{fele-def}\\   g_i^e&=\textcolor{black}{
  \int_{\partial K^e\cap\Gamma}\varphi_i\left(\mathbf{f}_h-
  \mathbf{f}(u_h)\right)\cdot\mathbf{n}\ds. }\label{gele-def}
 \end{align}
 
The time derivatives $\dot u_i$ of the Bernstein coefficients
corresponding to the standard
Galerkin approximation \eqref{galerkin1}
are given by the solution of the linear system
\beq
\label{galerkin2}
\sum_{j\in\mathcal N_i}m_{ij}\dot u_j=b_i(u_h,u_{\rm in})-
\sum_{e\in\mathcal E_i}\int_{K^e}
\varphi_i\nabla\cdot\mathbf{f}(u_h)\dx,
\qquad i=1,\ldots,N_h.
\eeq

By definition \eqref{prec} of the local mass lumping operator $P^e$, we have
$$\tilde C^e_k-C^e_k=P^eC^e_k-C^e_k=(M^e_L-M^e_C)(M^e_C)^{-1}C_k^e.$$
Using the global matrix/vector notation, the vector
$f^e=\{f_i^e\}_{i=1}^{N_h}$ of antidiffusive
element contributions defined by
\eqref{fele-def} can be written as
\beq\label{fele-def2}
f^e=(M_L^e-M_C^e)(\dot u+(M^e_C)^{-1}\mathbf{C}^e\cdot\mathbf{f})
-\textcolor{black}{\tilde{D}^e}u
-(\mathbf{C}^e)^\top\cdot\mathbf{f}+r^e,
\eeq
where $r^e$ is an element vector containing the contributions \textcolor{black}{
  $$r_i^e=\int_{K^e}\nabla\varphi_i\cdot\mathbf{f}(u_h)\dx
+\int_{\partial K^e\cap\Gamma_-}\varphi_i(u_h-u_i)
\mathbf{f}'(u_h)\cdot\mathbf{n}\ds.$$}

For any element vector $v^e\in\mathbb{R}^{N_h}$, the components of the matrix-vector products
$(M^e_L-M^e_C)v^e$ and $\black{\tilde{D}^e}v^e$ sum to zero. Moreover, the partition of unity
property of the Bernstein basis functions $\varphi_i$ implies that
$\sum_{i=1}^{N_h}\nabla\varphi_i=\mathbf{0}$
and, therefore, $\sum_{i=1}^{N_h}\mathbf{c}_{ji}^e=\mathbf{0}$ by definition 
\eqref{cij-def}. It follows that
\beq
\sum_{i=1}^{N_h}f_i^e=\sum_{i\in\mathcal N^e}f_i^e=
0\qquad\forall e=1,\ldots, E_h.
\eeq
 The full element matrices $M_C^e$ and
$\mathbf{C}^e$ can be calculated just once on
the reference element and multiplied by element-dependent
Jacobian data. A formula for
$\tilde{\mathbf{C}}^e$ is presented in the Appendix.
Note that the involved integrals $\sum_{\textcolor{black}{j\in\mathcal N^e\backslash\{i\}}}
  m_{ij}^e(\dot u_i-\dot u_j)=\int_{K^e}\varphi_i
(\dot u_i-\dot u_h)\dx$,
  $\sum_{j\in\mathcal N^e}\mathbf{c}_{ij}^e\cdot\mathbf{f}_j=
\int_{K^e}\varphi_i\nabla\cdot\mathbf{f}_h\dx$,
  and  $
  \sum_{j\in\mathcal N^e}\mathbf{c}_{ji}^e\cdot\mathbf{f}_j=
\int_{K^e}\nabla\varphi_i\cdot\mathbf{f}_h\dx$
can also be calculated directly in a matrix-free manner.
\medskip

In the next section, we decompose $f_i^e$ into a sum of
antidiffusive subcell fluxes $f_{ij}^e$ such that
\beq
f_i^e=\sum_{j\in\tilde{\mathcal N}_i^e} f_{ij}^{e},\qquad
 f_{ji}^{e}=-f_{ij}^{e}\quad\forall j\in\tilde{\mathcal N}_i^e.
 \eeq
Restricting the monolithic convex limiting strategy proposed in \cite{convex}
to $f_{ij}^{e}$, we will correct the bar states $\bar u_{ij}^e$ of the
low-order IDP scheme \eqref{bar} in a bound-preserving
manner. The limited counterpart
$f_{ij}^{e,*}$ of  $f_{ij}^{e}$
preserves the discrete conservation property and \red{is}
local extremum diminishing if
 \beq
 f_{ij}^{e,*}=0\quad\forall j\notin\tilde{\mathcal N}_i^e,\qquad
 f_{ji}^{e,*}=-f_{ij}^{e,*}\quad\forall j\in\tilde{\mathcal N}_i^e,
 \eeq
 \beq\label{ubarcorr}
 \bar u_{ij}^e\in\mathcal G\cap \mathcal G_i\quad\Rightarrow\quad
  \bar u_{ij}^{e,*}=\bar u_{ij}^e
 +\frac{f_{ij}^{e,*}}{2\black{\tilde{d}_{ij}^e}}\in\mathcal G\cap \mathcal G_i,
 \eeq
 where
 $\mathcal G_i$ is the set of states satisfying the
 local discrete maximum principle
 \beq\label{lmp-bar}
\min_{j\in\tilde{\mathcal N}_i}u_j=:u_i^{\min}\le\bar u\le u_i^{\max}:=
\max_{j\in\tilde{\mathcal N}_i}u_j.
\eeq
Note that we define the bounds $u_i^{\min}$ and $u_i^{\max}$ using the
subcell stencil
$\tilde{\mathcal N}_i=\bigcup_{e\in\mathcal E_i}\tilde{\mathcal N}_i^e
$
rather than the full
element stencil $\mathcal N_i$ of node $i$, \cyan{unless mentioned
  otherwise}. The pros and cons of \cyan{using tight bounds} are
explained in \cite{CG-BFCT} in the context of flux-corrected
transport (FCT) algorithms.
\smallskip

A locally bound-preserving IDP approximation to a given target flux $f_{ij}^e$
is given by \cite{convex}
\beq\label{fij-bar}
f_{ij}^{e,*}=\begin{cases}
  \min\,\left\{f_{ij}^e,\min\,\{2\black{\tilde{d}_{ij}^e}u_i^{\max}-\bar w_{ij}^e,
    \bar w_{ji}^e-2\black{\tilde{d}_{ij}^e}u_j^{\min}\}\right\} & \mbox{if}\ f_{ij}^e>0,\\[0.25cm]
    \max\left\{f_{ij}^e,\max\{2\black{\tilde{d}_{ij}^e}u_i^{\min}-\bar w_{ij}^e,
   \bar w_{ji}^e-2\black{\tilde{d}_{ij}^e}u_j^{\max}\}\right\} & \mbox{otherwise},
\end{cases}
\eeq
where \red{
$\bar w_{ij}^e=2\tilde{d}_{ij}^e\frac{u_j+u_i}{2}-\tilde{\bc}_{ij}^e
  \cdot(\bF_j-\bF_i)$. In infinite-precision arithmetic, this product has
  the same value as $2\tilde d_{ij}^e\bar u_{ij}^{e}$, where $\bar u_{ij}^e$
  is the bar state defined by \eqref{bar}. In numerical implementations,
  we calculate $\bar w_{ij}^e$ directly to avoid rounding errors due
  to division and multiplication by $\tilde{d}_{ij}^e$.}

\begin{rmk}
 Guermond and Popov \cite{Guermond2016} proved the validity of a local
 entropy inequality for \eqref{bar} using the fact that
 (see Theorem 4.7 in \cite{Guermond2016})
 $$
 E(\bar u_{ij}^e)\le\frac{E(u_i)+E(u_j)}{2}-\frac{\tilde{\mathbf{c}}_{ij}^e
   \cdot(\bfF(u_j)-\bfF(u_i))}{2\black{\tilde{d}_{ij}^e}}
 $$
 for any entropy pair $(E,\bfF)$.
   Our monolithic convex limiting strategy makes it possible to enforce such 
   inequality constraints for $E(\bar u_{ij}^{e,*})$ by reducing the magnitude
   of $f_{ij}^{e,*}$ if necessary. That is, the set $\mathcal G_i$
   may be redefined so as to enforce local entropy conditions in
   addition to local maximum principles.
\end{rmk}

\black{
\begin{rmk}
  Since the bar states $\bar u_{ij}$ of the low-order method for the
  linear advection equation with the flux function
  $\bff(\mathbf{x},u)=\bu(\mathbf{x}) u$
may fail to satisfy \eqref{maxbarstate}, the
  generalized version
\beq\label{fij-bar-linear}
f_{ij}^{e,*}=\begin{cases}
  \min\,\left\{f_{ij}^e,\max\{0,\min\,\{2\black{\tilde{d}_{ij}^e}u_i^{\max}-\bar w_{ij}^e,
    \bar w_{ji}^e-2\black{\tilde{d}_{ij}^e}u_j^{\min}\}\}\right\} & \mbox{if}\ f_{ij}^e>0,\\[0.25cm]
    \max\left\{f_{ij}^e,\min\{0,\max\{2\black{\tilde{d}_{ij}^e}u_i^{\min}-\bar w_{ij}^e,
   \bar w_{ji}^e-2\black{\tilde{d}_{ij}^e}u_j^{\max}\}\}\right\} & \mbox{otherwise}
\end{cases}
\eeq
of formula \eqref{fij-bar}
should be used to ensure positivity preservation for such linear flux functions.
\end{rmk}

To correct possible errors in the approximation of boundary terms, we define
\beq
b_i^*(u_h,u_{\rm in})=\tilde b_i(u_h,u_{\rm in})+\sum_{e\in\mathcal E_i}g_i^{e,*}
\eeq
using 
\beq
g_i^{e,*}=
\min\left\{g_i^{e,\max},\max\left\{g_i^e,g_i^{e,\min}\right\}\right\},
\eeq
where the target $g_i^e$ is defined by \eqref{gele-def} and the bounds are given by
\begin{align}
  g_i^{e,\max}&=(u_i^{\max}-u_i)
\int_{\partial K^e\cap\Gamma}
 \varphi_i|\mathbf{f}'(u_h)\cdot\mathbf{n}|\ds,\\
g_i^{e,\min}&=(u_i^{\min}-u_i)
\int_{\partial K^e\cap\Gamma}\varphi_i|\mathbf{f}'(u_h)\cdot\mathbf{n}|\ds.
\end{align}}
The semi-discrete version of the flux-corrected Galerkin
 scheme is given by
\beq\label{afc}
m_i\td{u_i}{t}=\sum_{e\in\mathcal E_i}
\sum_{j\in\tilde{\mathcal N}_i\backslash\{i\}}
[\black{\tilde{d}_{ij}^e}(u_j-u_i)+f_{ij}^{e,*}]-\sum_{j\in\tilde{\mathcal N}_i}
\tilde{\mathbf{c}}_{ij}\cdot\mathbf{f}_j+b_i^*(u_h,u_{\rm in}).
\eeq
The IDP property can be shown as before using the equivalent form
\begin{align*}
m_i\td{u_i}{t}
&=\sum_{e\in\mathcal E_i}
\sum_{j\in\tilde{\mathcal N}_i^e\backslash\{i\}}
2\black{\tilde{d}_{ij}^e}(\bar u_{ij}^{e,*}-u_i)+b_i^*(u_h,u_{\rm in}),\\
&=\sum_{e\in\mathcal E_i}
\sum_{j\in\tilde{\mathcal N}_i^e\backslash\{i\}}
2\black{\tilde{d}_{ij}^e}(\bar u_{ij}^{e,*}-u_i)+c_i(u_i^*-u_i)+\tilde b_i(u_h,u_{\rm in}),
\end{align*}
where $\bar u_{ij}^{e,*}$ is the flux-corrected bar state defined by
\eqref{ubarcorr}, $u_i^*\in\{u_i^{\min},u_i^{\max}\}$ and
\beq
0\le c_i\le \sum_{e\in\mathcal E_i}
\int_{\partial K^e\cap\Gamma}\varphi_i|\mathbf{f}'(u_h)\cdot\mathbf{n}|\ds
\eeq
by definition of
$g_i^{e,*}$. We remark that the representation of the flux-corrected scheme
in terms of $\bar u_{ij}^{e,*}$ and $u_i^*$ is used for
theoretical analysis only.
Practical implementations should be based on \eqref{afc}.

\begin{rmk}
In contrast to the element-based algorithms  
proposed in \cite{DG-BFCT,RD-BFCT,RD-BFCT2,CG-BFCT}, the above limiting strategy
rules out direct mass exchange between nodes that are not
nearest neighbors.
\end{rmk}


\begin{rmk}\label{remark:clipping_effects}
To avoid strong peak clipping effects and
achieve optimal convergence rates for 
$p>1$, the 
discrete maximum principle \eqref{lmp-bar} needs to be replaced
with less restrictive constraints in a neighborhood of smooth
local extrema \cite{CG-BFCT}. To that end,
a subset $\mathcal G_i$ of the
invariant set $\mathcal G$ can be defined, e.g.,
using the smoothness criteria presented in
\cite{Diot2012,dumbser2014,Guermond2018,hpfem,CG-BFCT}.
We explore this possibility further in \S\ref{sec:smoothness_indicator}.
\end{rmk}

\section{Computation of subcell antidiffusive fluxes}
\label{sec:subcell}

Clearly, the accuracy of the flux-corrected Galerkin discretization
\eqref{afc} depends on the
definition of the subcell fluxes $f_{ij}^e,\ j\in\tilde{\mathcal N}_i^e$
which we have left unspecified so far. The antidiffusive element
contributions defined by \eqref{fele-def} can be written as
\beq
f_i^e=
\sum_{j\in\tilde{\mathcal N}_i^e\backslash\{i\}}
\black{\tilde{d}_{ij}^e}(u_i-u_j)+q_i^e,
\eeq
where
 \begin{align}
q_i^e&=
\sum_{\textcolor{black}{j\in\mathcal N^e\backslash\{i\}}}
  m_{ij}^e(\dot u_i-\dot u_j)
  +\sum_{\textcolor{black}{j\in\mathcal N^e}}
  (\tilde{\mathbf{c}}_{ij}^e-\mathbf{c}_{ij}^e)\cdot\mathbf{f}_j
  -\sum_{\textcolor{black}{j\in \mathcal N^e}}\mathbf{c}_{ji}^e\cdot\mathbf{f}_j
\nonumber\\
&+
\int_{K^e}\nabla\varphi_i\cdot\mathbf{f}(u_h)\dx
\textcolor{black}{\ +
\int_{\partial K^e\cap\Gamma_-}\varphi_i(u_h-u_i)
   \mathbf{f}'(u_h)\cdot\mathbf{n}\ds}
 \label{fele-def3}
 \end{align}
 is the vector of  element contributions that require further decomposition into
 subcell fluxes.
\smallskip

The zero-sum property $\sum_{i=1}^{N_h}q_i^e=0$ of the element contributions
$q_i^e$ implies the existence of
a (generally non-unique) representation in the flux form
\beq
q_i^e=\sum_{j=1\atop i\ne j}^{N_h}q_{ij}^e,\qquad q_{ji}^e=-q_{ij}^e.\label{fdecomp}
\eeq
Let the auxiliary vector $v^e\in\mathbb{R}^N$ be defined as a solution of the linear system
\beq
(\hat M^{e}_L-\hat M^{e}_C)\hat v^e=\hat q^e,\label{fsolve}
\eeq
where $\hat q^e_{i_e}:=q_i^e$ for $i\in\mathcal N^e$.
The sparse $N\times N$ mass matrices
$$\hat M^{e}_C=\left\{\int_{K^e}\psi_i^e\psi_j^e\dx\right\}_{i,j=1}^N,
\qquad \hat M^{e}_L=\left\{\delta_{ij}\int_{K^e}\psi_i^e\dx\right\}_{i,j=1}^N
$$
are defined using the local basis functions $\psi_i^e$ of the piecewise
$\mathbb{P}_1/\mathbb{Q}_1$ B\'ezier net approximation on the macroelement
$K^e$. The subcell fluxes defined by
\beq
q_{ij}^e=\hat m_{i_ej_e}^{e}(\hat v_{i_e}^e-\hat v_{j_e}^e)
\eeq
satisfy \eqref{fdecomp} and vanish if nodes $i$ and $j$ are not nearest neighbors. The
matrix $\hat M^{e}_L-\hat M_C^{e}$  is symmetric with
vanishing row sums. Hence, the solution $\hat v^e$ of the auxiliary problem
\eqref{fsolve} is defined up to a constant. Since our definition of $q_{ij}^e$ is
independent of this constant, it can be chosen arbitrarily. In our
implementation, we solve \eqref{fsolve} subject to the linear
equality constraint $$\sum_{i=1}^N\hat v_i^e=0.$$

In summary, the original Galerkin discretization \eqref{galerkin1} can
be recovered using
\beq\label{ftarget}
f_{ij}^e=\black{\tilde{d}_{ij}^e}(u_i-u_j)+q_{ij}^e.
\eeq
 In
contrast to algebraic flux correction schemes for $\mathbb{P}_1$ and $\mathbb{Q}_1$
discretizations of general conservation laws \cite{Guermond2018,convex},
the error associated with the group
finite element approximation \eqref{groupf} cannot be neglected in high-order
versions. Our definition of the target fluxes $f_{ij}^e$ corrects this error
even for $p=1$.

\begin{rmk}
If the coefficients $\black{\tilde{d}_{ij}}$ of the graph Laplacian operator are defined
using the assembled global matrix $\tilde{\mathbf{C}}$,
the corresponding fluxes $f_{ij}$ should
be calculated using the formula
\beq
f_{ij}=\black{\tilde{d}_{ij}}(u_i-u_j)+\sum_{e\in\mathcal E_i}q_{ij}^e
\eeq
and limited using the low-order bar states
$\bar u_{ij}=\frac{u_j+u_i}{2}-\frac{\tilde{\bc}_{ij}
\cdot(\bF_j-\bF_i)}{2\black{\tilde{d}_{ij}}}$ of the global system.
\end{rmk}
  
\begin{rmk}
The 1D version of the compact-stencil FCT limiter introduced in
\cite{CG-BFCT}  is also based on a decomposition of generic element
contributions into (uniquely defined)
subcell fluxes. However, the multidimensional subcell
decomposition  proposed in Section 4.5 of \cite{CG-BFCT}
requires the computationally
intensive solution of
minimization problems and has not been tested in practice so far.
\end{rmk}

\section{Stabilization of subcell antidiffusive fluxes}\label{sec:stabilization}

The continuous Galerkin method exhibits suboptimal
$\mathcal O(h^p)$ convergence behavior even for smooth solutions of linear
advection problems on general meshes. To achieve optimal accuracy and
prevent formation of spurious ripples within the local bounds of the
limiting procedures, some high-order stabilization should be included
in the target flux. 
In the numerical studies of Lohmann et al.
\cite{CG-BFCT}, optimal convergence rates for high-order finite
element discretizations of the linear advection equation were
achieved using two-level Laplacian stabilization which can be
added to the vector $q^e$ before decomposing it into subcell
fluxes $q_{ij}^e$ in the manner described
in Section \ref{sec:subcell}. For nonlinear conservation laws, Guermond et al. \cite{Guermond2018,Guermond2014}
recommend the use of entropy viscosity (EV) stabilization. Its ability to
preserve the optimal order for $p>1$ is yet to be verified.
The same is true for stabilization via low-order approximations to the
nodal time derivatives $\dot u_i$, as proposed in \cite{convex}
for $p=1$.
    
The selection of genuinely high-order
stabilization tools for Bernstein finite element
approximations is beyond the scope of this work. In the
numerical experiments of \S\ref{sec:num}, we 
replace \eqref{ftarget} with
\beq\label{EV_flux}
  f_{ij}^{e,\rm stab}=(1-C_E\max\left(R_i,R_j\right))
  \tilde{d}_{ij}^{e}(u_i-u_j)+q_{ij}^e,
\eeq
where $R_i\in[0,1]$ is a nodal sensor that determines the appropriate
amount of nonlinear stabilization and $C_E=\mathcal{O}(1)$ is a user-defined
parameter (we use $C_E=1)$. 

Following Guermond et al. \cite{Guermond2018}, we choose an entropy pair
$(E(u),\bfF(u))$ for \eqref{ibvp-pde} and use
  \beq\label{entropy}
 R_i = \frac{\left|\sum_{j\in\mathcal N_i} [\bfF(u_j) - E^\prime(u_i)\bff(u_j)]\cdot \bc_{ij} \right|}
 {\left|\sum_{j\in\mathcal N_i} \bfF(u_j) \cdot \bc_{ij} \right| + \left|E^\prime(u_i)\right| \left|\sum_{j\in\mathcal N_i} \bff(u_j)\cdot \bc_{ij} \right| + \epsilon},
 \eeq
 where $\epsilon$ is a positive constant which prevents division by zero
 (we use $\epsilon=10^{-10}$). 
The so-defined $R_i$ measures the rate of entropy production
at node $i$. Note that we use
the coefficients $\bc_{ij}$ of the discrete gradient operator
corresponding to the high-order space in \eqref{entropy}.
This definition of $R_i$ extends the
domain of dependence to the full stencil $\mathcal N_i$
of node $i$ to improve robustness.
However, the stabilized subcell fluxes \eqref{EV_flux}
preserve the compact stencil property of the nonlinear
AFC scheme.

For all test problems in \S \ref{sec:num}, we use $E(u)=\frac{1}{2}u^2$
and $\bfF(u)=\int_0^u E^\prime(z)\bff^\prime(z)dz$. \red{For a detailed
discussion of entropy viscosity stabilization, we refer the
reader to Guermond et al. \cite{Guermond2018,Guermond2014}.}

\section{Extremum-preserving flux limiting}\label{sec:smoothness_indicator}
As mentioned in Remark \ref{remark:clipping_effects}, the local discrete
maximum principle \eqref{lmp-bar} may need to be relaxed to achieve high-order
convergence and alleviate peak clipping at smooth local extrema. In this
work, \red{we use one of the subcell smoothness indicators introduced
by Hajduk et al. \cite{RD-BFCT2}. The underlying smoothness criterion
is based on variations of the approximate nodal Laplacians
  \beq\label{laplacians}
  \tilde \eta_i=(\Delta_h\tilde u_h)_i
  :=\frac{1}{\tilde m_i}\int_\Omega \nabla\tilde u_h\cdot \nabla\psi_i\dx
  \eeq
  calculated using the piecewise $\mathbb{P}_1/\mathbb{Q}_1$
  basis functions $\psi_1,\ldots,\psi_{N_h}$, the diagonal entries 
  $\tilde m_i:=\int_\Omega \psi_i\dx$ of the corresponding lumped
  mass matrix, and the $\mathbb{P}_1/\mathbb{Q}_1$ interpolant 
  \beq
  \tilde u_h(\bfx)=\sum_{j=1}^{N_h}u_h(\mathbf{x}_j)\psi_j(\bfx)
  \eeq
  of $u_h(\mathbf{x}_j)=\sum_{k=1}^{N_h}\varphi_k(\bfx_j)u_k$, where
  $u_1,\ldots,u_{N_h}$ are
  the Bernstein degrees of freedom.
  Given the Laplacian reconstruction \eqref{laplacians},
  we calculate the nodal smoothness sensors \cite{RD-BFCT2}
  \begin{equation}\label{gamma1}
    \gamma_i = \begin{cases}
      \min\left\{1,
      \frac{C\max\{0,\eta_i^{\min}\eta_i^{\max}\}+\epsilon}{
      \max\{(\eta_i^{\min})^2,(\eta_i^{\max})^2\}+\epsilon}
      \right\} & \mbox{if } \mathbf x_i \in \Omega, \\
  1 & \mbox{if } \mathbf x_i \in \Gamma, \end{cases}
\end{equation}
  where $\epsilon>0$ is again a small positive number and $C\ge 1$
  is a sensitivity parameter. The maximum
$\eta_i^{\max}=\max_{j\in\tilde{\mathcal N}_i}\eta_j$ and minimum
$\eta_i^{\min}=\min_{j\in\tilde{\mathcal N}_i}\eta_j$
are taken over the set $\tilde{\mathcal N}_i$ of nodes
that share a subcell with node $i$.
Formula \eqref{gamma1} produces $\gamma_i=0$ if
the signs of $\eta_i^{\max}$ and $\eta_i^{\min}$
  differ. The maximal value $\gamma_i=1$ is attained if the signs of the
  two extremal values are the same and their magnitudes do not
  differ by more than a factor of $C$. In the numerical
studies below, we use $C=3$. 

To prevent unnecessary flux limiting at smooth peaks, we modify formula
\eqref{fij-bar-linear} as follows:
\beq\label{SI_flux}
f_{ij}^{e,*}=\begin{cases}
\min\,\{f_{ij}^e,\min\{\gamma_if_{ij}^{e}+(1-\gamma_i)
\max\{0,2\tilde{d}_{ij}^eu_i^{\max}-\bar w_{ij}^e\},\\
\phantom{\min\,\{f_{ij}^e,\min\{}
  \gamma_jf_{ij}^{e}+(1-\gamma_j)\max\{0,
  \bar w_{ji}^e-2\tilde{d}_{ij}^eu_j^{\min}\}\}\} & \mbox{if}\ f_{ij}^e>0,\\
\max\,\{f_{ij}^e,\max\{\gamma_if_{ij}^{e}+(1-\gamma_i)
\min\{0,2\tilde{d}_{ij}^eu_i^{\min}-\bar w_{ij}^e\},\\
\phantom{\max\,\{f_{ij}^e,\min\{}
  \gamma_jf_{ij}^{e}+(1-\gamma_j)\min\{0,
  \bar w_{ji}^e-2\tilde{d}_{ij}^eu_j^{\max}
  \}\}\}  & \mbox{otherwise}.
\end{cases}
\eeq
This modification relaxes the bounds of the flux constraints
associated with nodes $i$ and $j$ using the corresponding nodal
smoothness indicators. The IDP property w.r.t. the invariant
set $\mathcal G=[u^{\min},u^{\max}]$ can be enforced
by using the relaxed bounds $\gamma_iu^{\max}+(1-\gamma_i)u_i^{\max}$
and $\gamma_iu^{\min}+(1-\gamma_i)u_i^{\min}$ in the limiting 
formula \eqref{fij-bar-linear} instead of replacing it with
\eqref{SI_flux}, see \cite{RD-BFCT2} for details.
}

\section{Numerical examples}
\label{sec:num}

In this section, we apply the subcell flux limiting procedure to
(stabilized) Galerkin discretizations of scalar test problems. The
main purpose of this numerical study is to show that the proposed
low-order scheme and subcell flux decomposition are well suited
for algebraic flux correction purposes. More detailed studies
of stabilization approaches and smoothness indicators will
be presented elsewhere.

All computations are performed using {\sc Proteus} (https://proteustoolkit.org), an open-source Python toolkit for numerical simulations. 
We consider the following low-order methods:
\begin{itemize}
\item LO\,$\{$full stencil$\}$. In this version, we do not apply the mass
  lumping operator $P^e$ to the element matrices $C^e_k$ of the
  discrete gradient operator for  Bernstein elements of degree
  $p=1,2$. The  element matrix $\black{\tilde{D}^e}$ of the resulting discrete
  diffusion operator has $N^2$
  nonvanishing entries.
\item LO\,$\{$compact stencil$\}$.
  This is the low-order method defined by
  \eqref{nobar}. In this section, it is used for $p=2$ only. The
  element matrix $\tilde{D}^e$ of the discrete diffusion operator
  has the compact sparsity pattern of the piecewise $\mathbb{Q}_1$
  discretization on the 4-element submesh depicted in
  Fig.~\ref{fig:subcells}.
\end{itemize}
The high-order methods under investigation are abbreviated as follows:
\begin{itemize}
\item HO\,$\{$Galerkin, L$\}$.
  No stabilization of the Galerkin target \eqref{fele-def}, 
  limiting via \eqref{fij-bar} or \eqref{fij-bar-linear} for the nonlinear
  and the linear problems, respectively.
\item HO\,$\{$EV$\}$. Stabilized EV target \eqref{EV_flux}, no limiting.
\item HO\,$\{$EV,\,L$\}$. Stabilized EV target \eqref{EV_flux},
  limiting via \eqref{fij-bar} or \eqref{fij-bar-linear} for the nonlinear
  and the linear problems, respectively.
\item HO\,$\{$EV,\,L,\,SI$\}$. Stabilized EV target \eqref{EV_flux}, 
 limiting using the smoothness indicator \eqref{gamma1}.
\end{itemize}

In the rest of this section, we proceed as follows. We first consider 
linear advection problems which we solve using the full and compact stencil
versions of LO, as well as different versions of HO.
The objective  is to assess the quality of the low-order
method and to study the convergence behavior of the high-order method
in situations when the exact solution is smooth.
Thereafter, we solve two nonlinear problems using 
LO\,$\{$compact stencil$\}$, \mbox{HO\,$\{$Galerkin, L$\}$}, and
HO\,$\{$EV, L$\}$. The results of these numerical experiments
illustrate the IDP property
of the low-order method and the importance of using
high-order stabilization for the target fluxes.

\subsection{Linear advection}

\subsubsection{One-dimensional advection}
The first linear problem that we consider in this study
is the one-dimensional advection equation
\beq
\pd{u}{t}+v\pd{u}{x}=0\quad\mbox{in}\quad\Omega=(0,1)
\eeq
with the constant velocity $v=1$.
The smooth initial condition is given by 
\beq
u_0(x)=\exp\left[-100(x-0.25)^2\right].
\eeq
We solve this problem up to the final time $T=0.5$ and measure numerical
errors w.r.t. the $L^1$~norm.

The grid convergence history for the low-order methods under investigation
are reported in
Table~\ref{Table:oneD_lowOrder}. The experimental
orders of convergence (EOC) for pairs of uniform 1D meshes
are calculated using the formula presented in \cite{CG-BFCT}.
We observe that the accuracy of
the full stencil version deteriorates significantly
as we switch from the subcell
$\mathbb{Q}_1$ discretization to the $\mathbb{Q}_2$ approximation
with the same number of DoFs. The compact-stencil $\mathbb{Q}_2$
scheme produces more accurate results than its full-stencil counterpart.
The numerical studies presented in \cite{CG-BFCT} indicate that
more dramatic improvements can be expected for high-order
Bernstein elements. At least for constant velocities, the
convergence behavior of the compact-stencil version is largely
independent of $p$, as shown in \cite{CG-BFCT}.

In Table \ref{Table:oneD_highOrder}, we present the results of grid
convergence studies for the high-order stabilized $\mathbb{Q}_2$
approximations. \red{In the limited versions of the HO\,$\{$EV$\}$
  method, we use the
 full stencil bounds $u_i^{\max}=\max_{j\in\mathcal N_i}u_j$
and  $u_i^{\min}=\min_{j\in\mathcal N_i}u_j$.
It can be seen that the SI relaxation based
on \eqref{gamma1} and \eqref{SI_flux} results in smaller
global $L^1$ errors and faster convergence on coarse meshes.
However, the  EOCs of flux-limited approximations
are not as high as those of HO\,$\{$EV$\}$ in this 
example.}

\begin{table}[h]\scriptsize
  \begin{center}
    \begin{tabular}{|c||c|c||c|c||c|c||} \cline{1-7}
      \multicolumn{1}{|c||}{} &
      \multicolumn{2}{|c||}{LO\,$\{$full stencil$\}$, $\Q$} &
      \multicolumn{2}{|c||}{LO\,$\{$full stencil$\}$, $\QQ$} &
      \multicolumn{2}{|c||}{LO\,$\{$compact stencil$\}$, $\QQ$} \\ \hline
      $N_h$ &
      $\|u_h-u_\text{exact}\|_{L^1}$ & EOC &
      $\|u_h-u_\text{exact}\|_{L^1}$ & EOC  &
      $\|u_h-u_\text{exact}\|_{L^1}$ & EOC \\ \hline
      11  & 1.45E-2 & --   & 1.61E-2 & -- & 1.51E-2 & --   \\ \hline
      15  & 1.25E-2 & 0.44 & 1.42E-2 & 0.39 & 1.31E-2 & 0.42 \\ \hline
      20  & 1.07E-2 & 0.49 & 1.24E-2 & 0.44 & 1.13E-2 & 0.47 \\ \hline
      28  & 8.79E-3 & 0.56 & 1.03E-2 & 0.50 & 9.34E-3 & 0.54 \\ \hline
      39  & 7.10E-3 & 0.62 & 8.51E-3 & 0.57 & 7.59E-3 & 0.60 \\ \hline
      54  & 5.65E-3 & 0.68 & 6.89E-3 & 0.63 & 6.08E-3 & 0.66 \\ \hline
      75  & 4.40E-3 & 0.74 & 5.46E-3 & 0.69 & 4.77E-3 & 0.72 \\ \hline
      105 & 3.36E-3 & 0.79 & 4.23E-3 & 0.75 & 3.66E-3 & 0.78 \\ \hline
      147 & 2.52E-3 & 0.84 & 3.22E-3 & 0.80 & 2.76E-3 & 0.82 \\ \hline      
    \end{tabular}
    \caption{Linear advection in 1D,
      grid convergence history for the low-order methods.         
      \label{Table:oneD_lowOrder}}
    \bigskip
    
    \begin{tabular}{|c||c|c||c|c||c|c||} \cline{1-7}
      \multicolumn{1}{|c||}{} &
      \multicolumn{2}{|c||}{HO\,$\{$EV$\}$, $\QQ$} &
      \multicolumn{2}{|c||}{HO\,$\{$EV,L$\}$, $\QQ$} &
      \multicolumn{2}{|c||}{HO\,$\{$EV,L,SI$\}$, $\QQ$} \\ \hline
      $N_h$ &
      $\|u_h-u_\text{exact}\|_{L^1}$ & EOC &
      $\|u_h-u_\text{exact}\|_{L^1}$ & EOC &
      $\|u_h-u_\text{exact}\|_{L^1}$ & EOC \\ \hline
      11  & 5.09E-3 & --   & 7.77E-3 & --   & 6.50E-3 & --   \\ \hline
      15  & 3.05E-3 & 1.51 & 4.84E-3 & 1.40 & 3.63E-3 & 1.73 \\ \hline
      20  & 1.69E-3 & 1.94 & 2.66E-3 & 1.96 & 1.76E-3 & 2.37 \\ \hline
      28  & 7.39E-4 & 2.35 & 1.20E-3 & 2.27 & 7.89E-4 & 2.27 \\ \hline
      39  & 2.99E-4 & 2.64 & 6.33E-4 & 1.85 & 3.43E-4 & 2.43 \\ \hline
      54  & 1.25E-4 & 2.62 & 3.20E-4 & 2.05 & 1.47E-4 & 2.55 \\ \hline
      75  & 4.97E-5 & 2.76 & 1.54E-4 & 2.18 & 6.24E-5 & 2.56 \\ \hline
      105 & 1.87E-5 & 2.87 & 7.04E-5 & 2.30 & 2.63E-5 & 2.54 \\ \hline
      147 & 6.96E-6 & 2.91 & 3.35E-5 & 2.18 & 1.14E-5 & 2.46 \\ \hline
    \end{tabular}
    \caption{Linear advection in 1D,
      grid convergence history for the high-order methods.         
      \label{Table:oneD_highOrder}}
  \end{center}
\end{table}
\begin{rmk}
  To avoid errors due to inaccurate initialization, we $L^2$-project the smooth initial
  data of this test problem into the $\QQ$ finite element space by solving
  a linear system with the consistent mass matrix.  For all other test
  problems, we define the  Bernstein coefficients $u_i(0)=u_0(\mathbf{x}_i)$
  using the (generally inaccurate
  but bound-preserving, cf. \cite{phillipsInterp}) interpolation
  at the control points $\mathbf{x}_i$.
\end{rmk}

\subsubsection{Solid body rotation}\label{sec:num:sbr}
To facilitate a direct
comparison with the $\mathbb{P}_1/\mathbb{Q}_1$ version of
algebraic flux correction schemes and variational approaches
to shock capturing, let us now
consider the solid body rotation benchmark \cite{john2008,convex,afc1,leveque}.
In this 2D experiment, we solve the unsteady linear advection equation
  \begin{equation*}
    \pd{u}{t}+\nabla\cdot(\mathbf{v}u)=0\quad\mbox{in}\quad\Omega=(0,1)^2
  \end{equation*}
  using the divergence-free velocity field ${\bf v}(x,y)=2\pi(0.5-y,x-0.5)^\top$ to rotate a
  slotted cylinder, a sharp cone, and a smooth hump around the center
  $(0.5,0.5)$ of the domain $\Omega$.
  Homogeneous boundary conditions are prescribed on $\Gamma_-$.
  The initial condition, as defined by LeVeque \cite{leveque},
  is given by
  \begin{equation*}
    u_0(x,y)=\begin{cases}
    u_0^{\rm hump}(x,y)
    &\text{if}\  \sqrt{(x - 0.25)^2 + (y - 0.5)^2}\le 0.15, \\
    u_0^{\rm cone}(x,y)
    &\text{if}\ \sqrt{(x - 0.5)^2 + (y - 0.25)^2}\le 0.15, \\
    1 &\text{if}\ \begin{cases}
    \left(\sqrt{(x - 0.5)^2 + (y - 0.75)^2}\le 0.15 \right), \\
    \left(|x - 0.5| \ge 0.025,~ y\ge 0.85\right),
    \end{cases}\\
    0 &  \text{otherwise},
    \end{cases}
  \end{equation*}
  where\vspace{-0.25cm}
  \begin{align*}
    u_0^{\rm hump}(x,y)&=    \frac14 + \frac14 \cos \left(
    \frac{\pi \sqrt{(x - 0.25)^2 + (y - 0.5)^2}}{0.15}\right),\\
    u_0^{\rm cone}(x,y)&= 1-\frac{\sqrt{(x - 0.5)^2 + (y - 0.25)^2}}{0.15}.
  \end{align*}
  After each complete revolution, the exact solution coincides with the
  initial condition.

  \begin{figure}[h!]
    \centering
    \begin{tabular}{ccc}
      {\thead{$E_1=9.74\times 10^{-2}$ \\ $u^{\max}=0.5423$}} &
      {\thead{$E_1=1.05\times 10^{-1}$ \\ $u^{\max}=0.4730$}} &
      {\thead{$E_1=9.62\times 10^{-2}$ \\ $u^{\max}=0.5631$}} \\
      \includegraphics[scale=0.2]{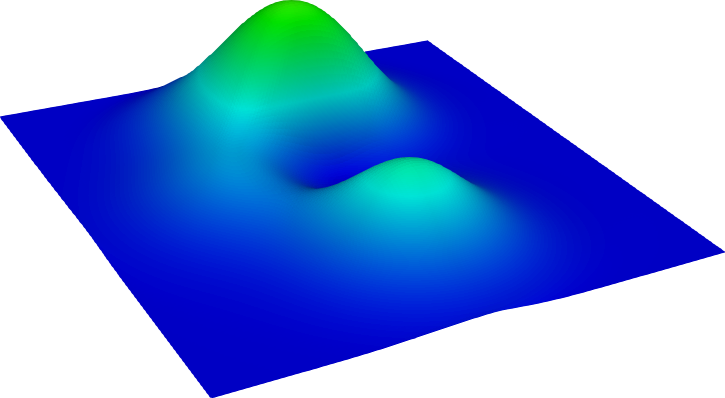} &
      \includegraphics[scale=0.2]{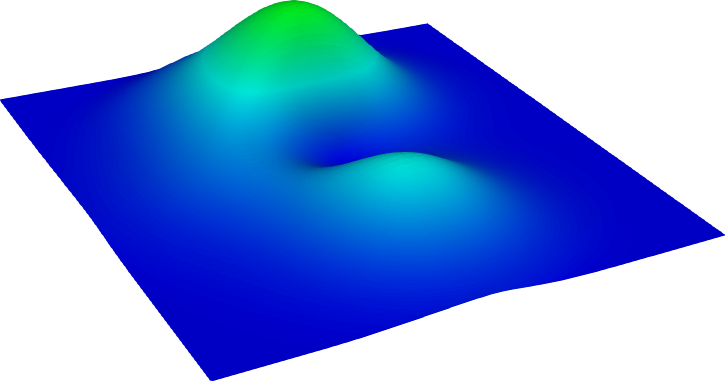} &
      \includegraphics[scale=0.2]{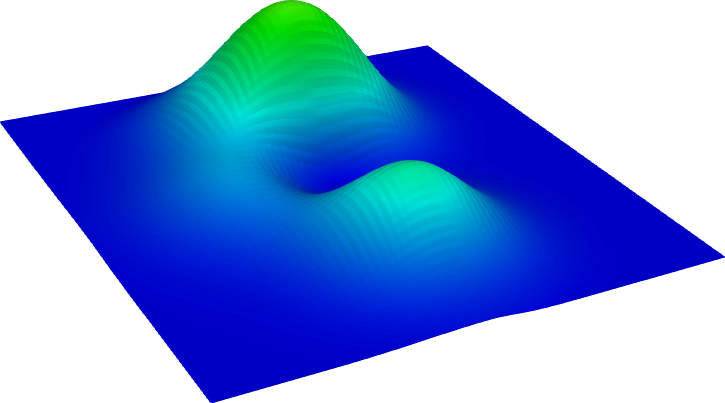} \\
      && \\
      {\thead{$E_1=7.92\times 10^{-2}$ \\ $u^{\max}=0.6650$}} &
      {\thead{$E_1=8.80\times 10^{-2}$ \\ $u^{\max}=0.6220$}} &
      {\thead{$E_1=7.80\times 10^{-2}$ \\ $u^{\max}=0.6631$}} \\
      \includegraphics[scale=0.2]{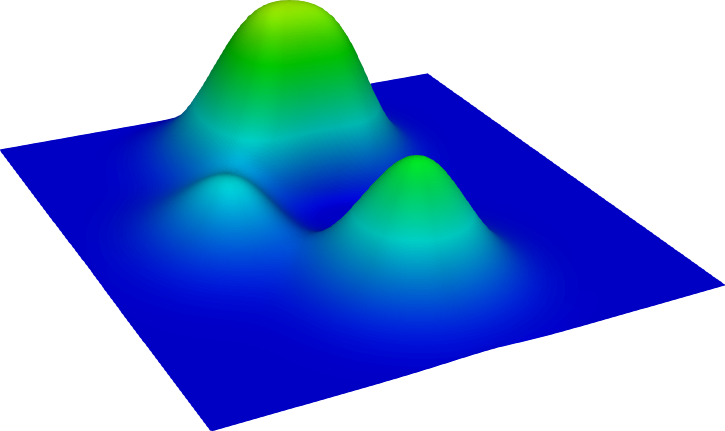} &
      \includegraphics[scale=0.2]{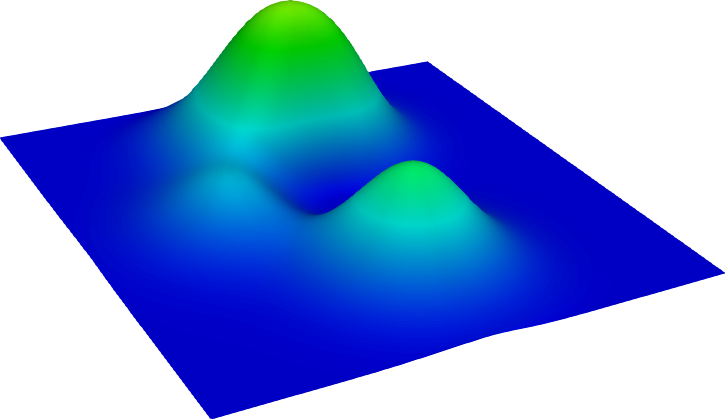} &
      \includegraphics[scale=0.2]{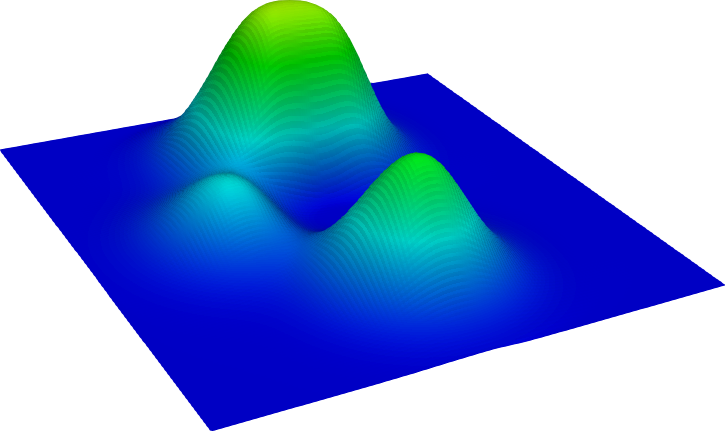} \\
      && \\
      LO\,$\{$full stencil$\}$, $\Q$ &
      LO\,$\{$full stencil$\}$, $\QQ$ &
      LO\,$\{$compact stencil$\}$, $\QQ$ 
    \end{tabular}
    \caption{Solid body rotation \cite{leveque}. Low-order solutions
      after one full rotation (final time $T=1$). The total number
      of DoFs is $N_h=129^2$ in the diagrams of the first row and
      $N_h=257^2$ in the diagrams of the second row.
      \label{fig:solid_rotation_low_order}}
  \end{figure}

  In Figure \ref{fig:solid_rotation_low_order}, we show the low-order
 $\Q$ and $\QQ$ approximations at the final time $T=1$ (one full
  rotation). The diagrams of the first and second row were obtained
  using $N_h=129^2$ and $N_h=257^2$ DoFs, respectively.
  For a better quantitative comparison,
  the $L^1$ errors $E_1=\|u_h-u_\text{exact}\|_{L^1}$ and
  the global maxima $u^{\max}=\max_{i=1,\dots, N_h} u_i$ of
  the Bernstein coefficients
  are listed above each plot. As expected, the approximation calculated
  with the full stencil $\QQ$ scheme proves more dissipative than
  the compact-stencil $\Q$ and $\QQ$ approximations. In contrast
  to the subcell upwinding strategy employed in
  \cite{RD-BFCT,CG-BFCT}, the low-order scheme defined
  by \eqref{nobar} preserves the  $\Q$ sparsity pattern exactly
  even for nonuniform velocity fields and nonlinear flux functions.
  This remarkable property eliminates a major bottleneck to
  achieving high performance and $p$-independent convergence
  behavior with matrix-based algebraic flux correction schemes.
  In our numerical 
  experiment, the low-order $\QQ$ solution obtained with
  \eqref{nobar} is as accurate as the subcell $\Q$
  approximation with the same number of DoFs.

  \begin{figure}[h!]
    \centering
    \begin{tabular}{cccc}
      {\thead{$E_1=2.01\times 10^{-2}$ \\ $u^{\max}=0.9868$}} &
      {\thead{$E_1=3.21\times 10^{-2}$ \\ $u^{\max}=1.0158$}} &
      {\thead{$E_1=3.46\times 10^{-2}$ \\ $u^{\max}=0.9562$}} &
      {\thead{$E_1=3.32\times 10^{-2}$ \\ $u^{\max}=0.9752$}} \\
      \includegraphics[scale=0.15]{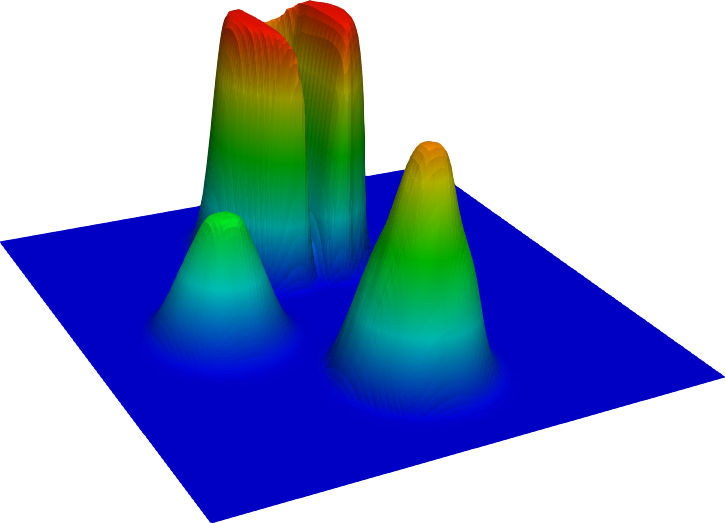} &
      \includegraphics[scale=0.15]{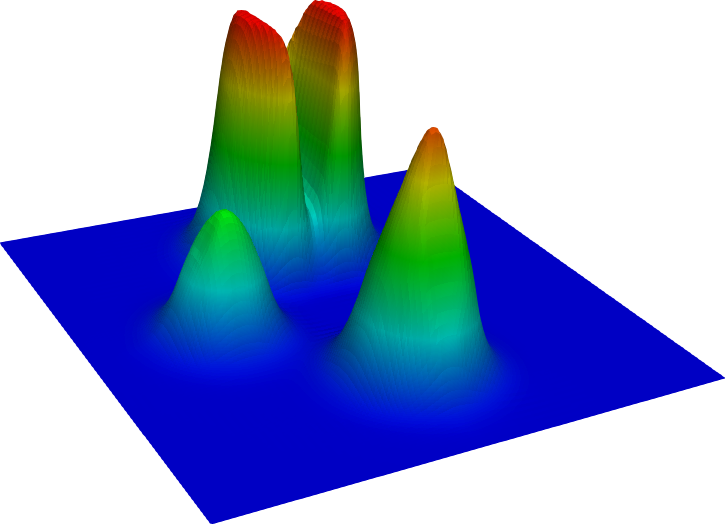} &
      \includegraphics[scale=0.15]{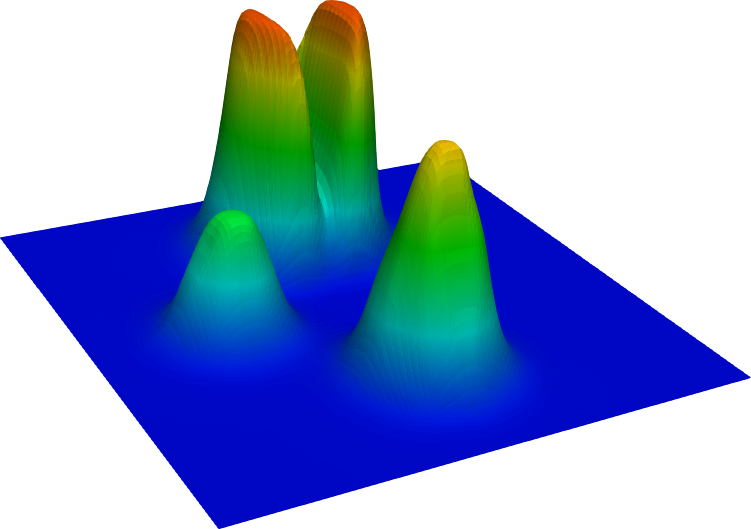} &
      \includegraphics[scale=0.15]{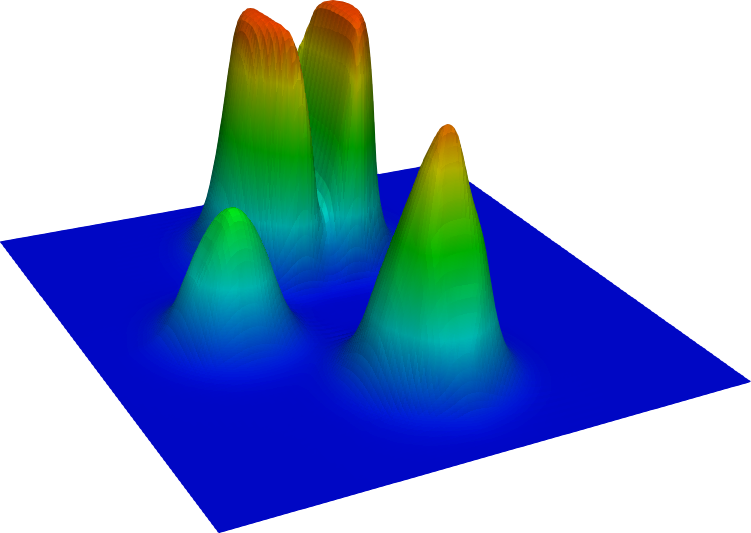} \\
      &&& \\
      {\thead{$E_1=1.12\times 10^{-2}$ \\ $u^{\max}=0.9996$}} &
      {\thead{$E_1=1.87\times 10^{-2}$ \\ $u^{\max}=1.0126$}} &
      {\thead{$E_1=2.03\times 10^{-2}$ \\ $u^{\max}=0.9945$}} &
      {\thead{$E_1=1.94\times 10^{-2}$ \\ $u^{\max}=0.9964$}} \\
      \includegraphics[scale=0.15]{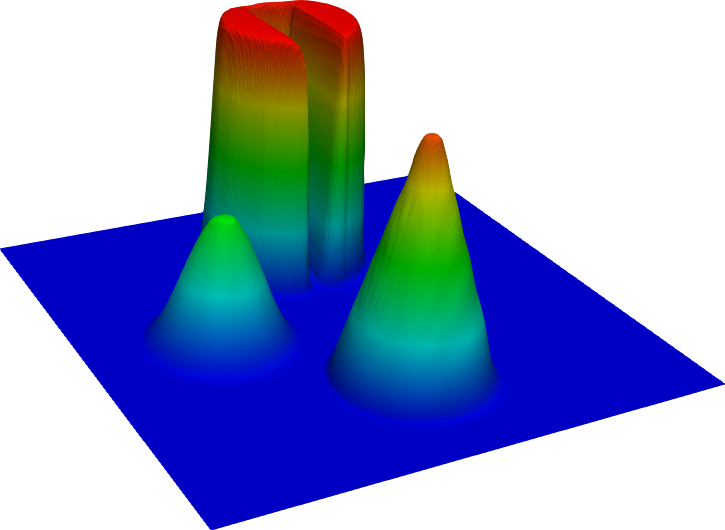} &
      \includegraphics[scale=0.15]{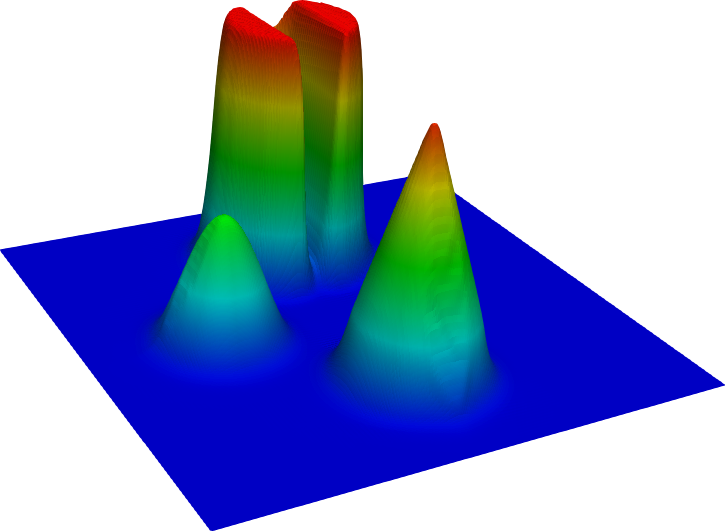} &
      \includegraphics[scale=0.15]{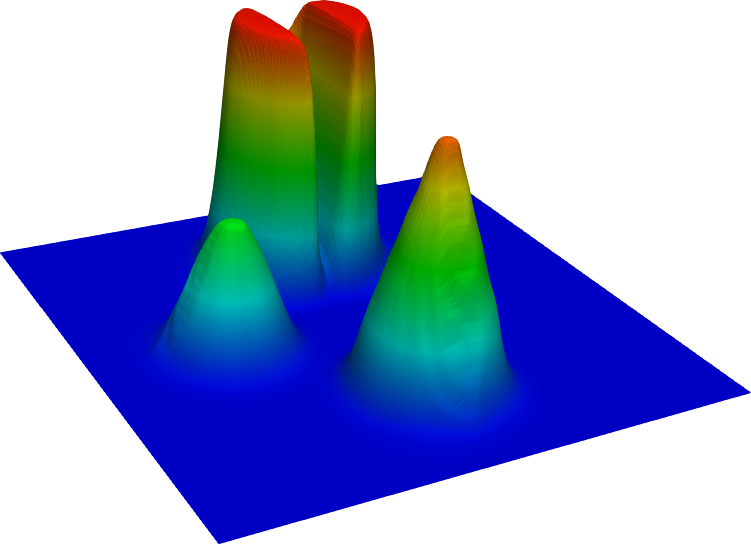} &
      \includegraphics[scale=0.15]{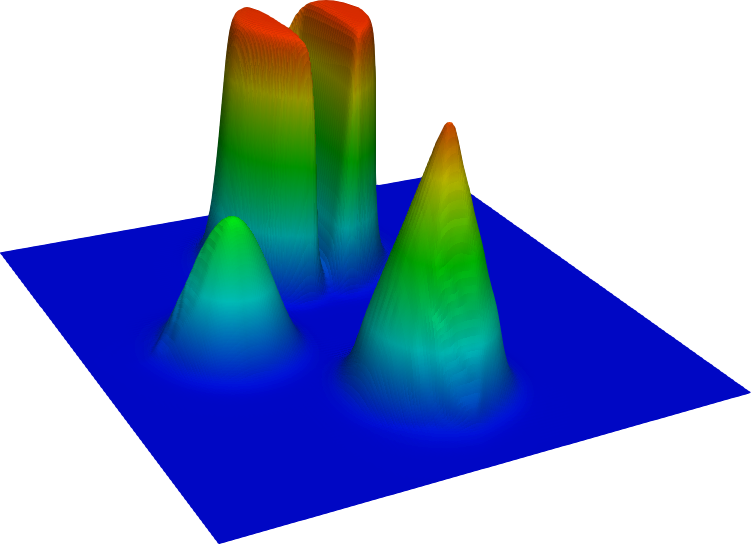} \\
      &&& \\
     HO\,$\{$Galerkin,L$\}$ &
     HO\,$\{$EV$\}$ &
     HO\,$\{$EV,L$\}$ &
     HO\,$\{$EV,L,SI$\}$
    \end{tabular}
    \caption{Solid body rotation problem \cite{leveque}.
     High-order solutions
      after one full rotation (final time $T=1$). The total number
      of DoFs is $N_h=129^2$ in the diagrams of the first row and
      $N_h=257^2$ in the diagrams of the second row.
      \label{fig:solid_rotation_high_order}}
  \end{figure}
    
  \begin{figure}[!h]
    \centering
    \subfloat[$N_h=129^2$]{
      \begin{tabular}{cc}
      \includegraphics[scale=0.07]{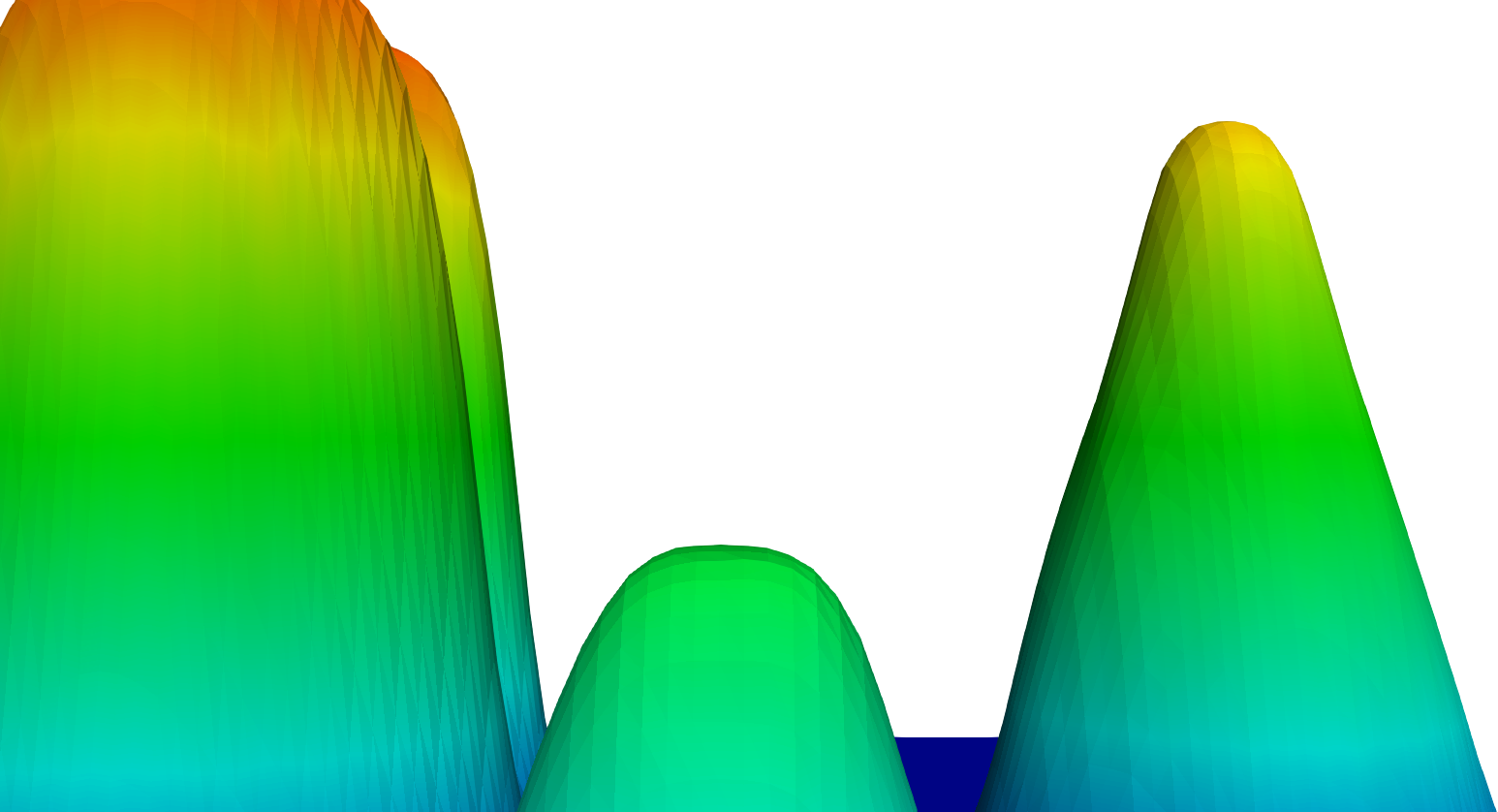} &
      \includegraphics[scale=0.07]{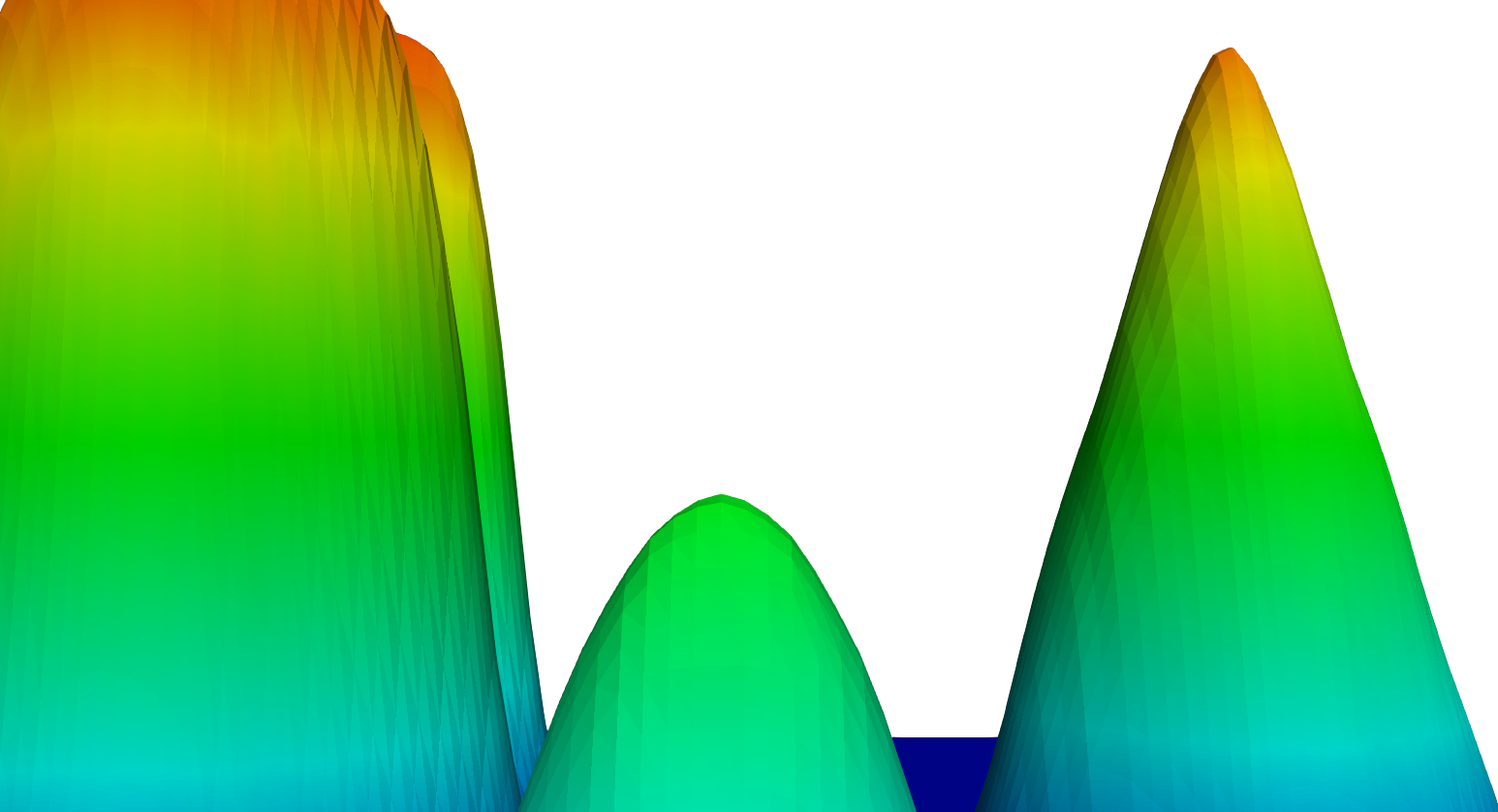} \\
      & \\
      HO\,$\{$EV,L$\}$ &
      HO\,$\{$EV,L,SI$\}$
      \end{tabular}
    }
    ~
    \subfloat[$N_h=257^2$]{
      \begin{tabular}{cc}
      \includegraphics[scale=0.07]{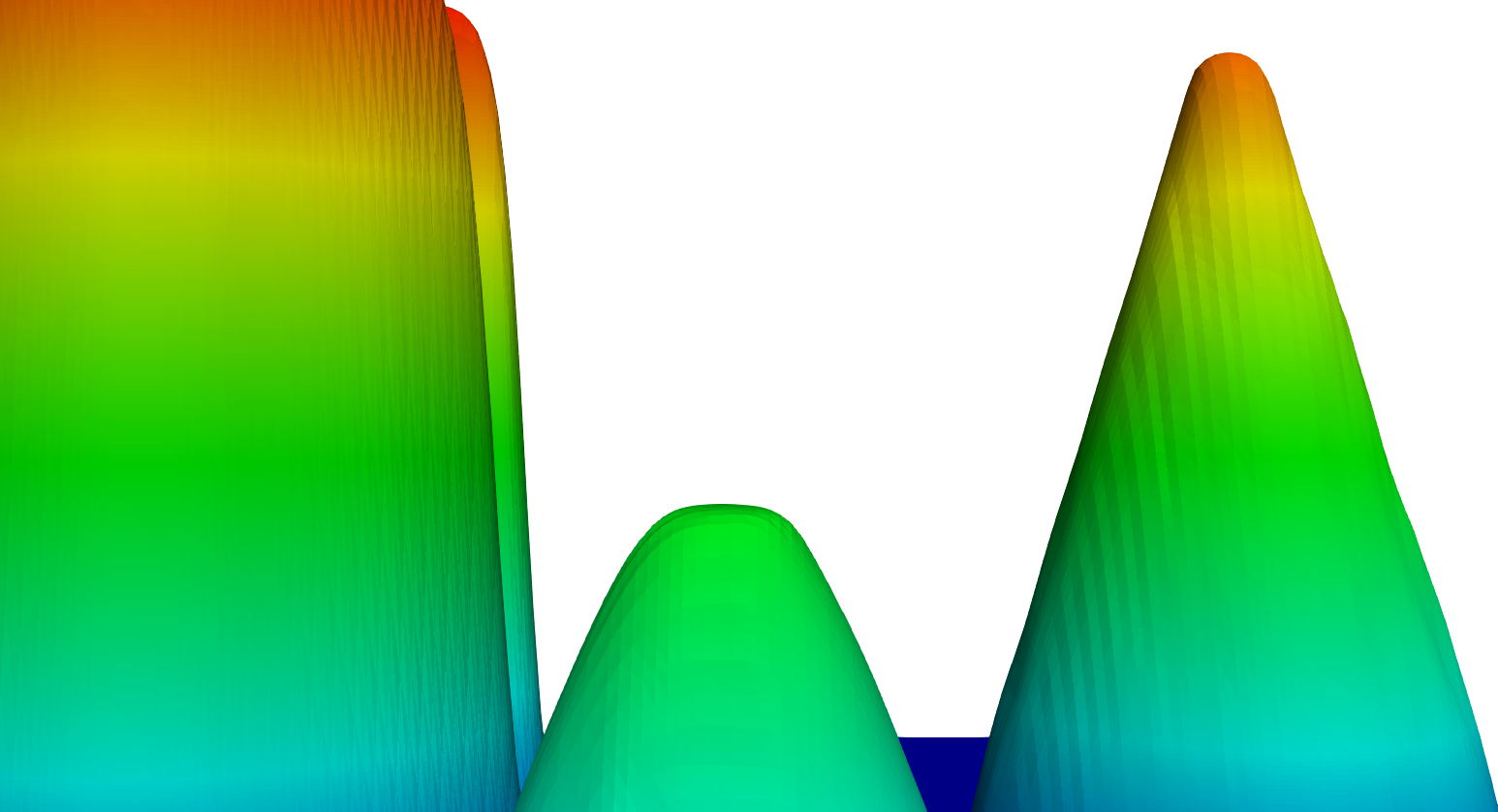} &
      \includegraphics[scale=0.07]{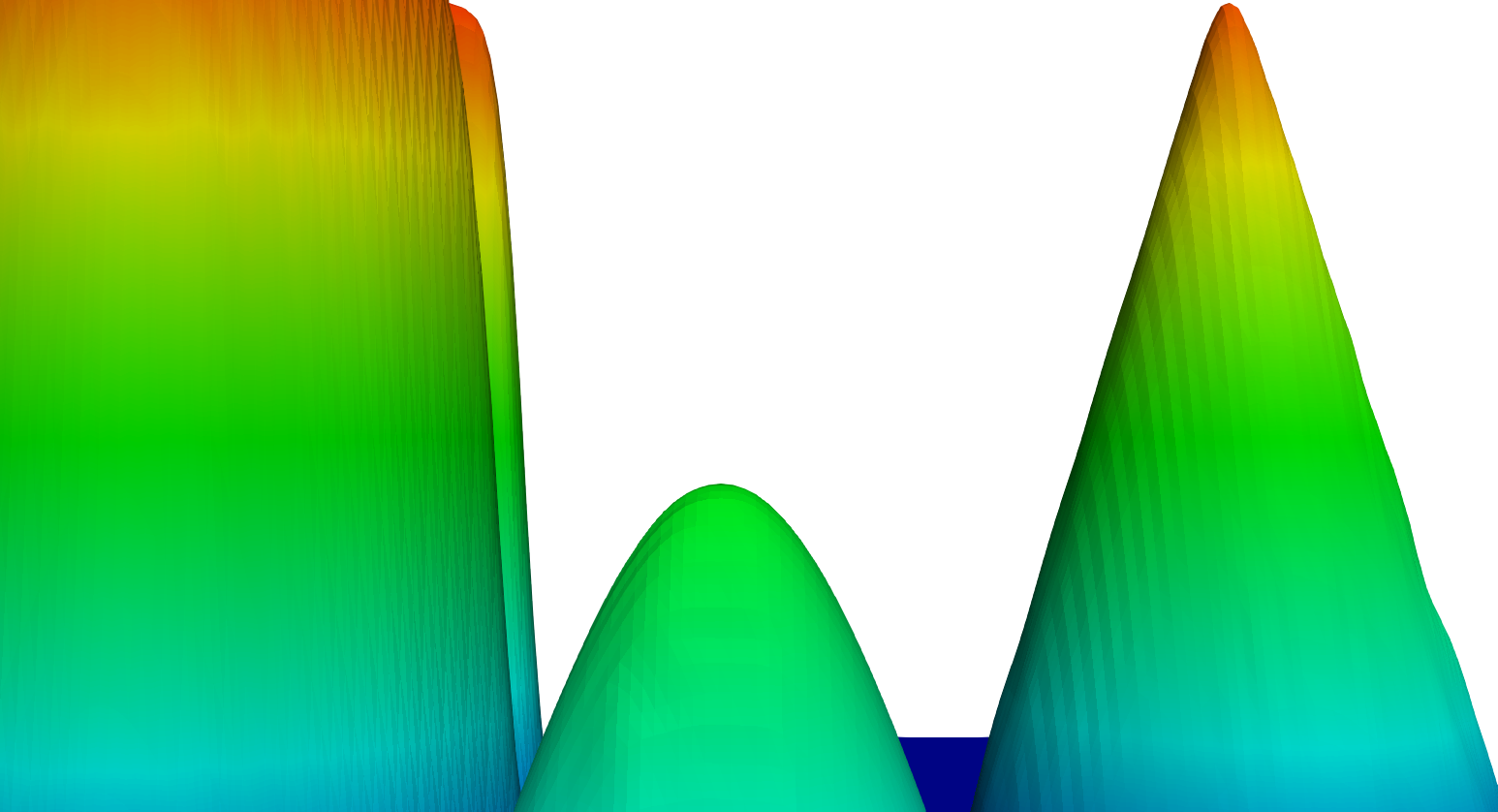} \\
      & \\
      HO\,$\{$EV,L$\}$ &
      HO\,$\{$EV,L,SI$\}$
      \end{tabular}
      }
    \caption{Solid body rotation problem \cite{leveque}. Zooms of the
     limited high-order $\QQ$ solutions at $T=1$
     obtained without and with using the smoothness indicator
     defined by \eqref{gamma1} to
     reduce peak clipping effects.
      \label{fig:solid_rotation_high_order_zoomed}}
  \end{figure}

  The results obtained with the high-order extensions and zooms
  of the flux-limited solutions are shown in Figs
  \ref{fig:solid_rotation_high_order} and
  \ref{fig:solid_rotation_high_order_zoomed}, respectively.
  The  activation of subcell flux correction eliminates small undershoots
  and overshoots at the edges of the slotted cylinder but smears
  the bound-preserving peaks of the hump and cone significantly.
  The stabilization of the target flux via entropy viscosity
  increases the $L^1$ error without having any positive impact
  on the quality of the flux-corrected $\QQ$ approximations in
  this particular example.
  \red{The Laplacian-based smoothness indicator
    $\gamma_i$ defined by \eqref{gamma1} was used to relax the
  bounds in formula \eqref{SI_flux}.}
  The SI version recognizes the top of the hump as a smooth
  extremum and resolves it very well even on the coarser mesh.
  Flux limiting at the top of the cone is deactivated as soon
  as the peak becomes rounded enough for \eqref{gamma1} to produce
  \red{$\gamma_i^e=1$}. At the same time, no violation of discrete
  maximum principles occurs in the neighborhood of 
  discontinuities, where the second derivatives exhibit abrupt
  changes and \eqref{gamma1} produces \red{$\gamma_i^e=0$.}

    \subsubsection{Steady circular advection}
\begin{figure}[h]
    \centering
    \begin{tabular}{ccc}
      {\thead{$E_1=9.14\times 10^{-2}$ \\ $u_h\in[0,1]$}}&
      {\thead{$E_1=1.75\times 10^{-2}$ \\ $u_h\in [0,1]$}} &
      {\thead{$E_1=5.61\times 10^{-2}$ \\ $u_h\in [0,1]$}} \\
      \includegraphics[scale=0.2]{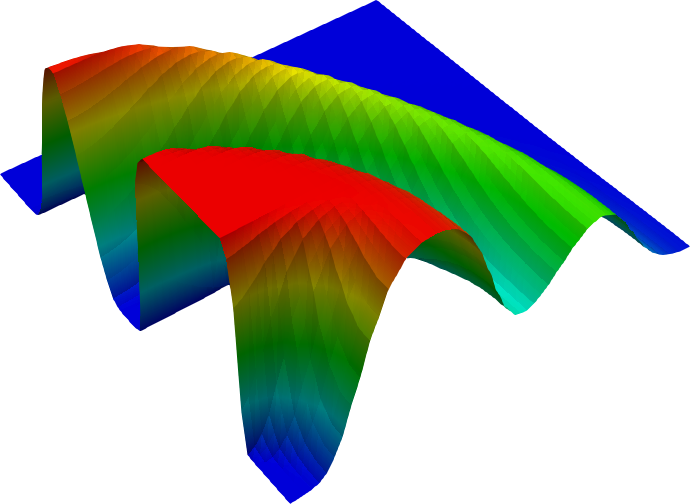} &
      \includegraphics[scale=0.2]{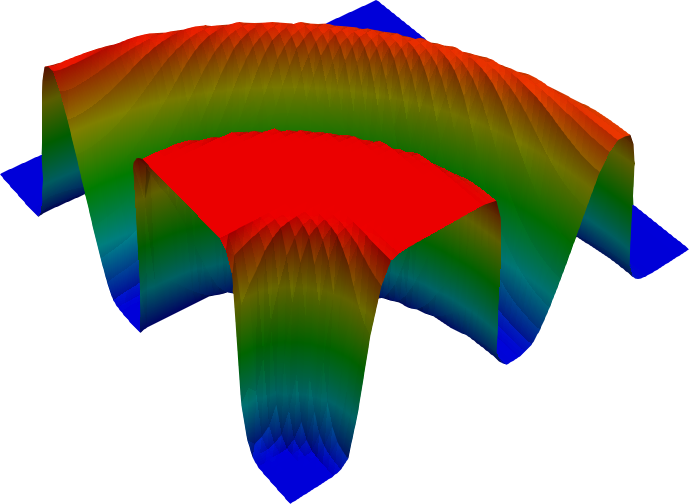} &
      \includegraphics[scale=0.2]{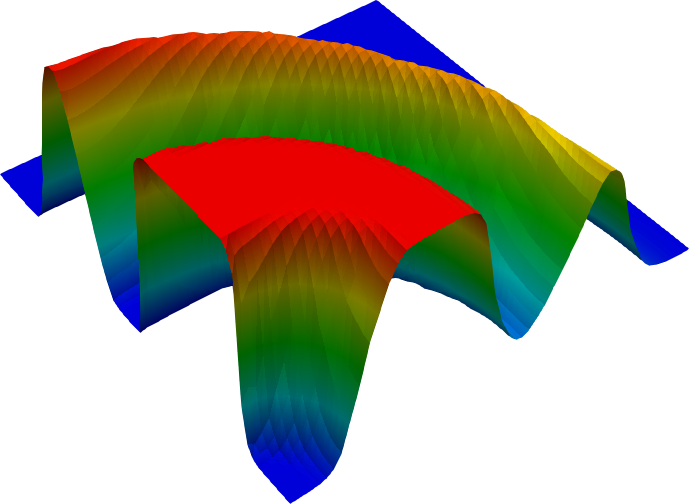} \\
      && \\
      {\thead{$E_1=5.84\times 10^{-2}$ \\ $u_h\in [0,1]$}}&
      {\thead{$E_1=7.96\times 10^{-3}$ \\ $u_h\in [0,1]$}} &
      {\thead{$E_1=3.03\times 10^{-2}$ \\ $u_h\in [0,1]$}} \\
      \includegraphics[scale=0.2]{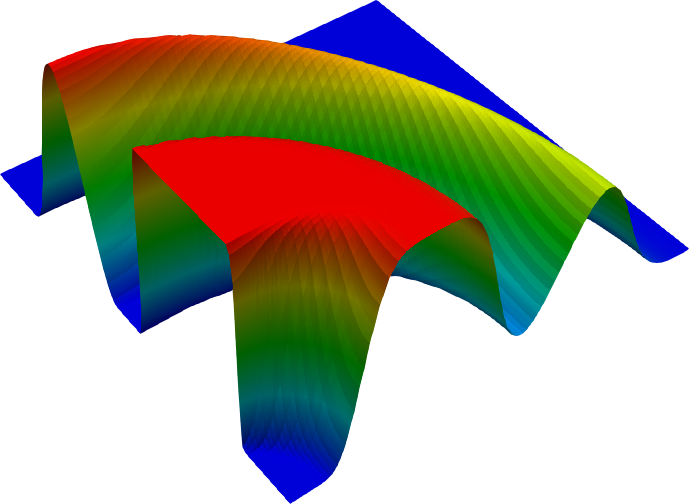} &
      \includegraphics[scale=0.2]{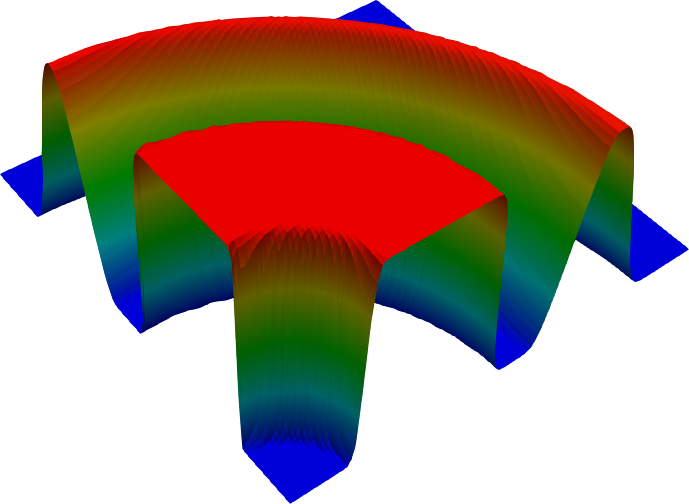} &
      \includegraphics[scale=0.2]{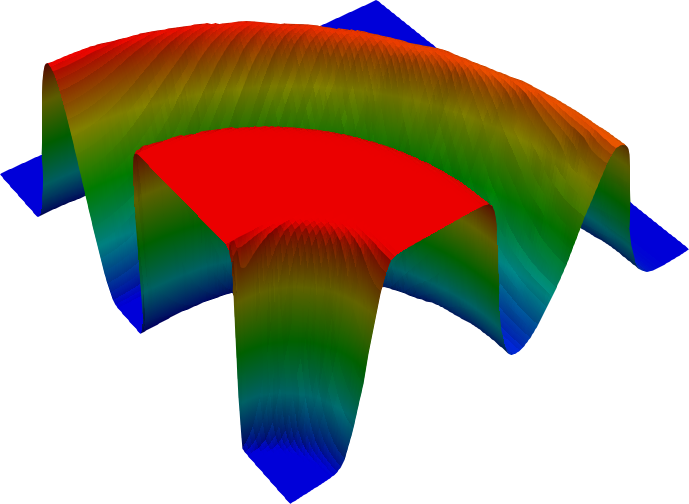} \\      
      && \\
      LO\,$\{$compact stencil$\}$ &
      HO\,$\{$Galerkin,L$\}$ &
      HO\,$\{$EV,L$\}$
    \end{tabular}
    \caption{Steady circular advection. Stationary $\QQ$ solutions
      calculated using time marching. The total number of DoFs is
      $N_h=65^2$ in the diagrams of the first row and
    $N_h=129^2$ in the diagrams of the second row.
      \label{fig:steady_state}}
  \end{figure}

    In contrast to the FCT algorithms employed in
    \cite{DG-BFCT,Guermond2018,Guermond2019,RD-BFCT,CG-BFCT},
    the monolithic convex limiting strategy is well suited for
    calculating steady-state solutions. To show this,
     we solve 
\begin{equation}
\nabla\cdot(\mathbf{v}u)=0\quad\mbox{in}\  
\Omega=(0,1)^2\label{cc_eq}
\end{equation}
using the divergence-free velocity field
$\mathbf{v}(x,y)=(y,-x)$. The inflow
boundary condition and the exact solution at any point
in $\bar\Omega$ are given by
\begin{equation}
  u(x,y)=\left\{\begin{array}{ll}
  1, &\quad \mbox{if} \ \ 0.15\le r(x,y)\le 0.45,\\
  \cos^2\left(10\pi\frac{r(x,y)-0.7}{3}\right),
  &\quad \mbox{if} \ \ 0.55\le r(x,y)\le 0.85,\\
  0, & \quad\mbox{otherwise},
  \end{array}\right. \label{cc-init}
\end{equation}
where $r(x,y)=\sqrt{x^2+y^2}$ denotes the distance to the corner point $(0,0)$. 
The stationary $\QQ$ solutions obtained
with $N_h=65^2$ and $N_h=129^2$ are shown in Fig.
\ref{fig:steady_state}. These numerical solutions were marched to
the steady state by solving the lumped-mass version of the $\QQ$
approximation to the time-dependent advection problem
until the prescribed tolerance was reached
for the steady-state residuals.

\subsection{Burgers equation}

As a first nonlinear  test problem, we consider the
2D inviscid Burgers equation \cite{GuermondNazarov2014,convex}
\beq
\pd{u}{t}+\nabla\cdot\left(\mathbf{v}\frac{u^2}{2}\right)=0\qquad
\mbox{in}\ \Omega=(0,1)^2, 
\eeq
where $\mathbf{v}=(1,1)^\top$ is a constant vector. The
piecewise-constant initial data is given by
\beq
u_0(x,y)=\begin{cases}
-0.2 & \mbox{if}\quad x < 0.5\ \land y >0.5,\\
-1.0 & \mbox{if}\quad x > 0.5\ \land y >0.5,\\
\phantom{-}0.5 & \mbox{if}\quad x < 0.5\ \land y <0.5,\\
\phantom{-}0.8 & \mbox{if}\quad x > 0.5\ \land y <0.5.
\end{cases}
\eeq
The inflow boundary conditions are defined using the exact solution of the
pure initial value problem in~$\R^2$.  This solution can be found in
\cite{GuermondNazarov2014} and stays in the invariant set
$\mathcal G=[-1.0,0.8]$.

\begin{figure}[h!]
  \centering
  \begin{tabular}{ccc}
    ${\thead{E_1=1.94\times 10^{-2} \\ u_h\in [-1,0.8]}}$ & 
    ${\thead{E_1=1.16\times 10^{-2} \\ u_h\in [-1,0.8]}}$ &
    ${\thead{E_1=1.16\times 10^{-2} \\ u_h\in[-1,0.8]}}$
                \\
    \includegraphics[scale=0.25]{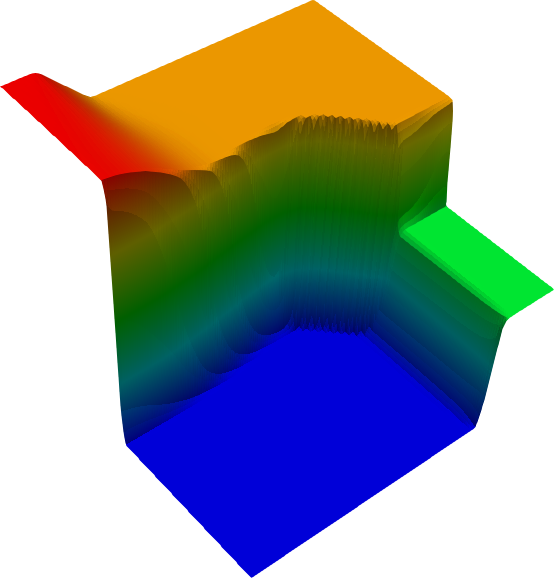}  &
    \includegraphics[scale=0.25]{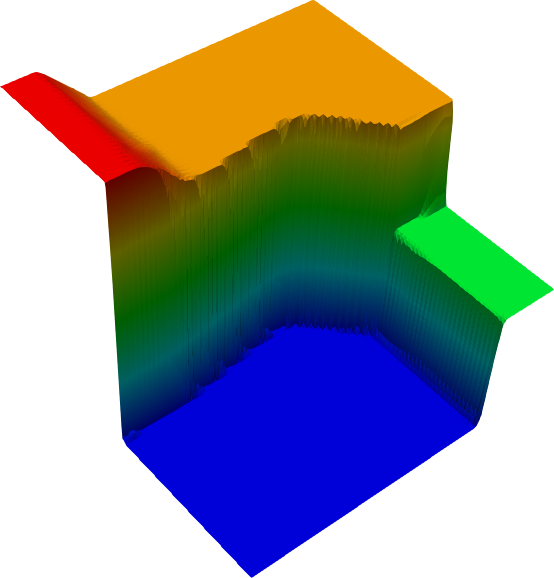} &
    \includegraphics[scale=0.25]{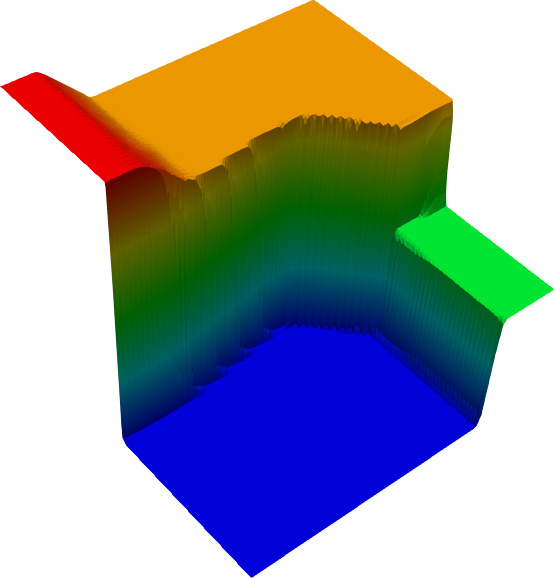}
    \\
    &&\\
    ${\thead{E_1=1.06\times 10^{-2} \\ u_h\in [-1,0.8]}}$ & 
    ${\thead{E_1=6.10\times 10^{-3} \\ u_h\in [-1,0.8]}}$ &
    ${\thead{E_1=6.10\times 10^{-3} \\ u_h\in [-1,0.8]}}$
    \\
    \includegraphics[scale=0.25]{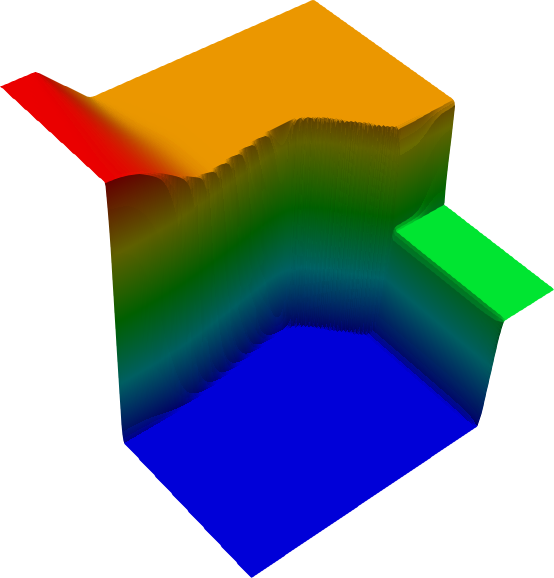}  &
    \includegraphics[scale=0.25]{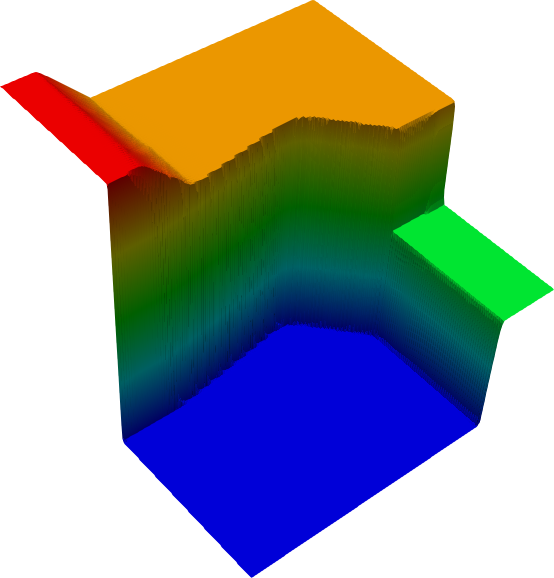} &
    \includegraphics[scale=0.25]{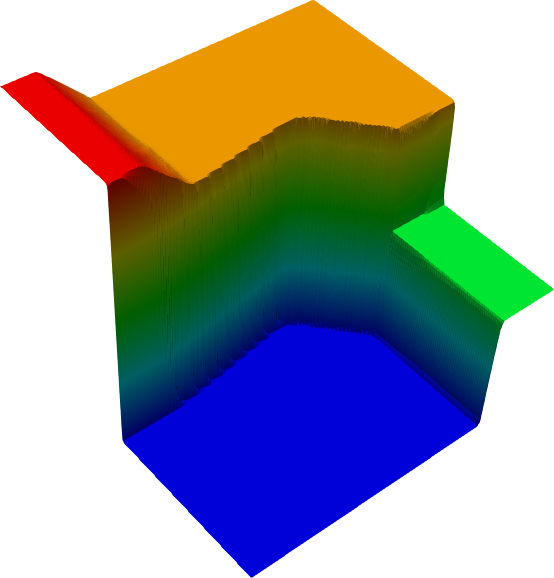}
    \\ 
    &&\\
    LO\,$\{$compact stencil$\}$ &
    HO\,$\{$Galerkin,L$\}$&
    HO\,$\{$EV,L$\}$ 
  \end{tabular}
  \caption{Burgers equation, bound-preserving $\QQ$ approximations at $T=0.5$.
The total number of DoFs is
      $N_h=129^2$ in the diagrams of the first row and
$N_h=257^2$ in the diagrams of the second row.
\label{fig:burgers}}
\end{figure}

The numerical solutions obtained at $T=0.5$ using $\QQ$ elements
with $N_h=129^2$ and $N_h=257^2$ DoFs are shown in Fig.~\ref{fig:burgers}.
The presented $L^1$ errors indicate that considerable amounts of
numerical diffusion can be safely removed in the process of subcell flux
correction. The use of EV stabilization has no significant
impact on the accuracy of the flux-corrected HO solutions in this example.

\subsection{KPP problem}

The KPP problem \cite{Guermond2014,Guermond2016,Guermond2017,kpp} is a
more challenging nonlinear test. In this final 2D experiment,
we solve the scalar conservation law
\eqref{ibvp-pde} with the nonlinear and
nonconvex flux function
\beq
\mathbf{f}(u)=(\sin(u),\cos(u))
\eeq
 in the domain
$\Omega=(-2,2)\times(-2.5,1.5)$ using the initial
condition
\beq
u_0(x,y)=\begin{cases}
\frac{14\pi}{4} & \mbox{if}\quad \sqrt{x^2+y^2}\le 1,\\
\frac{\pi}{4} & \mbox{otherwise}.
\end{cases}
\eeq
A simple (but rather pessimistic) upper bound for the
 maximum speed is $\lambda=1$. More
accurate GMS bounds can be found in \cite{Guermond2017}.
The exact solution exhibits a two-dimensional rotating
wave structure. The main challenge of this test is to
prevent possible convergence to wrong weak solutions.

\begin{figure}[!h]
  \centering
  \begin{tabular}{cccc}
    \includegraphics[scale=0.2]{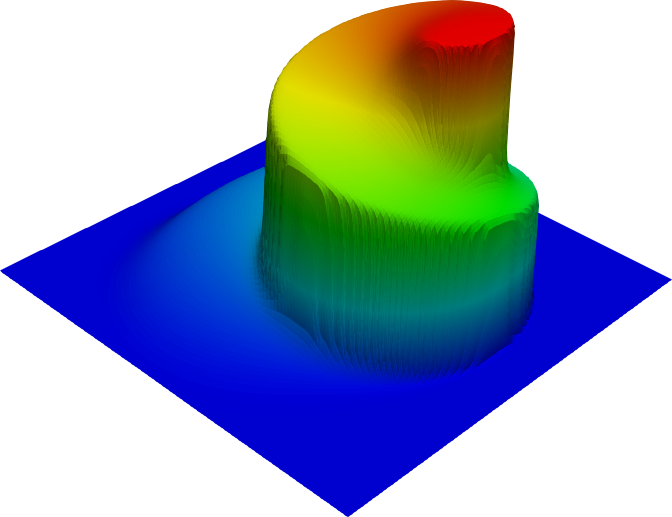} &
    \includegraphics[scale=0.2]{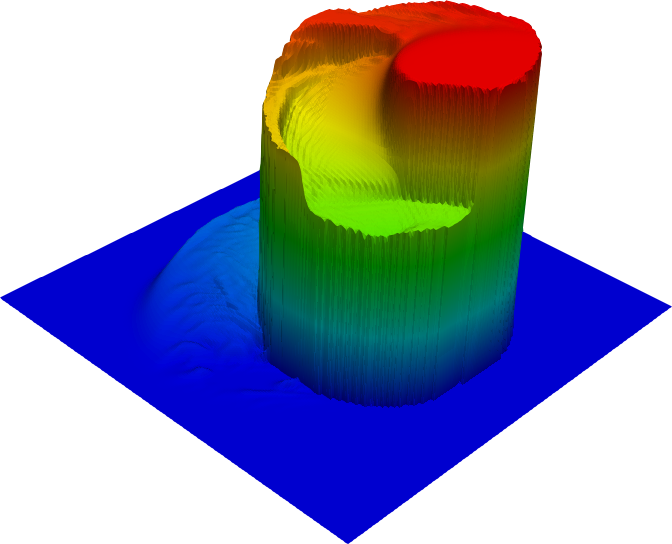} &
    \includegraphics[scale=0.2]{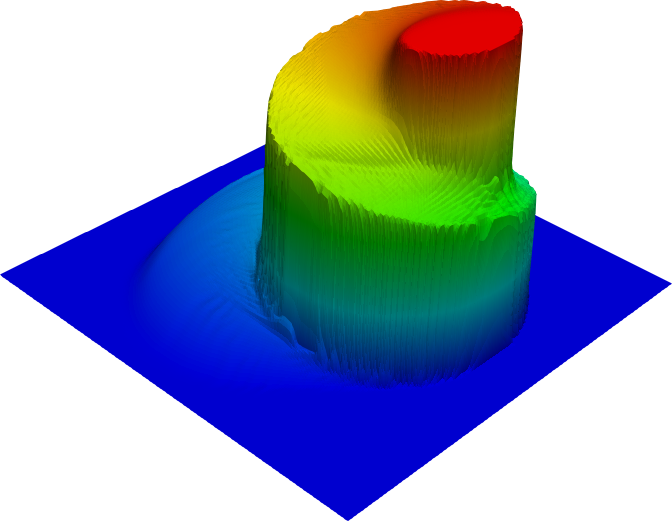}
    \\ 
    &&\\
     LO\,$\{$compact stencil$\}$ &
     HO\,$\{$Galerkin,L$\}$&
     HO\,$\{$EV,L$\}$ 
  \end{tabular}
  \caption{KPP problem \cite{kpp}, bound-preserving
    $\QQ$ approximations at $T=1$.
The total number of DoFs is $N_h=257^2$. 
  \label{fig:KPP}}
\end{figure}

The numerical solutions obtained at $T=1$ using $N_h=257^2$ DoFs are
displayed in Fig.~\ref{fig:KPP}. The plot shown in the middle
demonstrates that the lack of nonlinear stabilization in the
target flux of the AFC scheme may, indeed, cause convergence to an
entropy-violating solution. This example confirms the 
findings of Guermond et al. \cite{Guermond2014} who observed such
nonphysical behavior of flux-limited Galerkin methods in the
context of predictor-corrector FCT algorithms. The use of
entropy viscosity stabilization in the EV target 
of the monolithic AFC discretization
cures this drawback without introducing inordinately large
amounts of numerical dissipation (compare the well-resolved
solution on the right of Fig.~\ref{fig:KPP}
to the diffusive and
distorted approximations shown in the other two diagrams).

\section{Conclusions}\label{sec:conclusions}

The main result of this work is the development of a novel subcell flux correction
procedure for high-order finite elements. The proposed definitions of the low-order scheme
and of the antidiffusive fluxes lead to compact-stencil approximations which can be implemented
efficiently. The monolithic convex limiting strategy ensures well-posedness of nonlinear
discrete problems and opens new avenues for theoretical analysis of  high-order
AFC schemes. Since the $\mathbb{P}_1$ and $\mathbb{Q}_1$ versions of the presented
methodology have already been successfully applied to the
Euler equations of gas dynamics in \cite{convex}, it is hoped that
extensions of subcell flux limiting
to high-order Bernstein finite element discretizations
of nonlinear hyperbolic systems will be relatively straightforward.

\medskip
\paragraph{\bf Acknowledgments}

The work of Dmitri Kuzmin was supported by the German Research Association (DFG) under grant KU 1530/23-1.
\black{
  The work of Manuel Quezada de Luna was supported by King Abdullah University of Science
  and Technology (KAUST) in Thuwal, Saudi Arabia.}
The authors would like to thank Prof. David I. Ketcheson (KAUST) and
Christoph Lohmann (TU Dortmund University) for helpful discussions.


\section*{Appendix: Sparsity of the lumped discrete gradient}

Let $\lambda_1,\ldots,\lambda_{d+1}:\hat K\to [0,1]$ denote the
barycentric coordinates (i.e., $\mathbb{P}_1$ basis functions)
associated with the vertices of a $d$-dimensional reference
element $\hat K$.
The 1D Bernstein basis functions of degree $p\in\mathbb{N}$
are defined on the unit interval $\hat K=[0,1]^d$ thus:
$$
B_\alpha^p(x)=\left(\begin{array}{c}p\\ \alpha\end{array}\right)
  \lambda_1^\alpha(x)\lambda_2^{p-\alpha}(x),\qquad 0\le \alpha\le p.
$$
  The corresponding tensor
  product basis for  $\mathbb{Q}_p(\hat K)$ on a unit $d$-box
  $\hat K=[0,1]^d$ is defined as follows:
  $$
  B_\alpha^p(x_1,\ldots,x_d)=B_{\alpha_1}^p(x_1)\ldots B_{\alpha_d}^p(x_d),
  \qquad 0\le \alpha_1,\ldots,\alpha_d\le p.
  $$
  The Bernstein basis of $\mathbb{P}_p(K^e)$ 
  a $d$-simplex $K^e=\mathrm{conv}\{\mathbf{x}_1^e,
 \ldots,\mathbf{x}_{d+1}^e\}$ is given by
$$
B_\alpha^p(\lambda_1,\ldots,\lambda_{d+1})=
\frac{p!}{\alpha_1!\cdot\ldots\cdot \alpha_{d+1}!}
\lambda_1^{\alpha_1}\cdot\ldots\cdot\lambda_{d+1}^{\alpha_{d+1}},$$
where  $\alpha=(\alpha_1,\ldots,\alpha_{d+1})$ is a multiindex such that
$$
|\alpha|:=\alpha_1+\ldots+ \alpha_{d+1}=p.
$$

Consider the $N\times N$ element matrices $P^e=M_L^e(M_C^e)^{-1}$ and
$C_k^e,\ k=1,\ldots,d$ of the polynomial space
spanned by $\varphi_i^e=B^p_{\alpha(i)},\ i=1,\ldots,N.$
By definition \eqref{precC}, the $j$-th column of the
element matrix $\tilde C_k^e$ contains the Bernstein
coefficients of $\pd{\varphi_j^e}{x_k}$ multiplied by the diagonal entries
$m_i^e=\frac{|K^e|}{N}$ of $M_L^e$ \cite{CG-BFCT}. Indeed, the Bernstein polynomial
$B_h=\sum_{i=1}^N\tilde c_{ij,k}^e\varphi_i^e$ is the unique solution of
$$
\int_{K^e}\varphi_h^eB_h\hat{\dx}
=\frac{|K^e|}{N}
\int_{\hat K}\varphi_h^e\pd{\varphi_j^e}{x_k}\hat{\dx},\qquad
\varphi_h^e\in\{\varphi_1^e,\ldots,\varphi_N^e\}.$$
The solution of this linear system yields the Bernstein coefficients
of the local $L^2$ projection
which is exact for polynomials of degree up to
$p$. It follows that $B_h=\frac{|K^e|}{N}\pd{\varphi_j^e}{x_k}$.
\smallskip

Using the product rule, the gradient of the Bernstein basis function
$B_\alpha^p$ on a $d$-dimensional simplex element
 $\hat K$ can be written as
\cite{Kirby,Kirby2014}
$$
\nabla B_\alpha^p=p\sum_{k=1}^{d+1}B_{\alpha-e_k}^{p-1}\nabla\lambda_k.
$$
The degree elevation formula for simplicial Bernstein elements
yields
$$
B_{\alpha-e_k}^{p-1}=\frac{1}{p}
\sum_{l=1}^{d+1}(\alpha_k-e_k+e_l+1)B^p_\alpha
$$
and the compact sparsity pattern follows from the fact that
$$\nabla B_\alpha^p=
\sum_{k=1}^{d+1}
\sum_{l=1}^{d+1}
(\alpha_k-e_k+e_l+1)\nabla\lambda_kB_{\alpha-e_k+e_l}^p=
\sum_{|\beta|=p}\tilde c_\beta B_{\beta}^p,
$$
where $\tilde c_\beta=0$ if $\beta\ne \alpha-e_k+e_l$ for
some $k,l\in\{1,\ldots,d\}$. Hence, the coefficient $\tilde
c_{ij,k}^e$  is nonvanishing
only if $j=i$ or $i$ and $j$ are nearest neighbors belonging
to the same grid line of the B\'ezier net.
\bigskip

To verify the sparsity of the element matrix $\tilde C^e$ for a multidimensional
$d$-box $K^e$, note that
$$\pd{B_\alpha^p}{x_k}=\pd{B_{\alpha_k}^p}{x_k}\prod_{l=1\atop
  l\ne k}^d B_{\alpha_l}^p,\qquad k=1,\ldots, d.$$
The desired result follows from
the proof of sparsity for the one-dimensional simplex element.
\smallskip

We remark that the above formulas can also be used for practical
calculation of the lumped discrete gradient operator.
Efficient algorithms for calculating and inverting the element matrices
of high-order Bernstein finite element spaces can be found in
\cite{Ainsworth2011,Kirby2011,Kirby2014,Kirby2017}.

\end{document}